\documentclass[11pt]{amsart}
\usepackage{amssymb,latexsym,amsmath}
\usepackage{graphics}

\input{RBsymbols.mac}
\newcommand{\ord}{\text{ord}}
\newcommand{\Tr}{\text{Tr}}

\theoremstyle{plain}
\newtheorem{theorem}{Theorem}[section]
\newtheorem{lemma}[theorem]{Lemma}
\newtheorem{proposition}[theorem]{Proposition}

\theoremstyle{remark}
\newtheorem{definition}[theorem]{Definition}
\newtheorem{example}[theorem]{Example}
\newtheorem{remark}[theorem]{Remark}

\numberwithin{equation}{section}

\def\Z{{\mathbb Z}}

\def\R{{\mathbb R}}

\def\Qp{{\mathbb Q}_p}

 {\\ \begin{minipage}{10cm}\begin{eqnarray*}\label{#1}}%
 {\end{eqnarray*}\end{minipage}\addtocounter{equation}{1}%
  \hfill (\theequation) \bigskip \\}

 {\begin{eqnarray*}\label{#1}}%
 {\end{eqnarray*}\addtocounter{equation}{1}}

                      {\hspace*{\fill} \solidbox  \bigskip \\}

\begin{document}

\title{A wavelet theory for local fields and related groups}

\author{John J.~Benedetto}
\address{Department of Mathematics\\
           University of Maryland\\
           College Park, MD 20742, USA} 
\email{jjb@math.umd.edu}
\urladdr{http://www.math.umd.edu/\textasciitilde jjb}          
\author{Robert L.~Benedetto}
\address {Department of Mathematics and Computer Science\\
           Amherst College\\
           Amherst, MA 01002-5000}
\email{rlb@cs.amherst.edu}
\urladdr{http://www.cs.amherst.edu/\textasciitilde rlb}

\thanks
{The first named author gratefully acknowledges support from 
NSF DMS
Grant 0139759 and ONR Grant N000140210398.}

\thanks
{The second named author gratefully acknowledges support from 
NSF DMS
Grant 0071541.}

\keywords
{wavelet, locally compact abelian group, p-adic field}
\subjclass[2000]
{Primary: 42C40; Secondary: 11S99, 42C15, 43A70}

\date{September 9, 2003}

\begin{abstract}
Let $G$ be a locally compact abelian group with compact open
subgroup $H$. The best known example of such a group is $G={\mathbb Q}_p$, the
field of $p$-adic rational numbers (as a group under addition), which has
compact open subgroup $H={\mathbb Z}_p$, the ring of $p$-adic integers.
Classical wavelet theories, which require a non-trivial discrete subgroup
for translations, do not apply to $G$, which may not have such a subgroup.
A wavelet theory is developed on $G$
using coset representatives of the discrete quotient
$\widehat{G}/H^{\perp}$ to circumvent
this limitation.
Wavelet bases are constructed by means of an iterative method
giving rise to so-called wavelet sets in the dual group $\widehat{G}.$
Although the Haar and Shannon wavelets are naturally antipodal
in the Euclidean setting,
it is observed that their analogues for $G$ are equivalent.
\end{abstract}

\maketitle

\newcounter{bean}


\section{Introduction}\label{sec1}
\renewcommand{\theequation}{\thesection.\arabic{equation}}
\setcounter{equation}{0}

\subsection{Background}\label{sbsec1.1}

Every locally compact abelian group (LCAG) $G$ is topologically
and algebraically
isomorphic to ${\mathbb R}^d \times G_o,$ where ${\mathbb R}^d$ is Euclidean
space and $G_o$ is a LCAG containing a compact open subgroup $H_o;$ note
that a subgroup $H$ of $G$ is open if and only if
$G/H$ is discrete. Because of applications, e.g., in signal
processing, and because of the role of periodization, most wavelet theory
has been constructed on ${\mathbb R}^d$ or on LCAGs containing cocompact
discrete subgroups, e.g., see \cite{Mey90} (1990), \cite{Dau92} (1992), 
\cite{Dah94} (1994), \cite{Luk94} (1994), 
\cite{Far97} (1997),
\cite{Mal98} (1998)
for the multiresolution analysis (MRA) theory, see \cite{Zak96} (1996),
\cite{DL98} (1997),
\cite{DLS97} (1997), \cite{HWW97} (1997), \cite{SW98} (1998),
\cite{BL99} (1999),
\cite{BL01} (2001), \cite{BS02} (2002)
for the theory of wavelet sets, and see
\cite{PS98} (1998), \cite{BMM99} (1999), \cite{BM99} (1999),
\cite{Cou99} (1999),
\cite{Pap03} (2003) for unifying general approaches.
There has also been work on wavelet theory outside
of the Euclidean setting, e.g.,
\cite{Lem89} (1989),
\cite{Dah95} (1995),
\cite{Hol95} (1995),
\cite{Tri96} (1996),
\cite{Alt97} (1997),
\cite{Joh97} (1997),
\cite{AV99} (1999),
\cite{HLPS99} (1999),
\cite{ST99} (1999),
\cite{AAG00} (2000),
\cite{AKLT00} (2000),
and classical work on
harmonic analysis on local fields, e.g.,
\cite{Tai75} (1975) and \cite{DP83} (1983).
In this paper, we
shall deal with wavelet theory
for functions defined on groups $G_o.$ Our results are outlined in
Subsection~\ref{sbsec1.3}.

We are motivated by our interest in nondiscrete locally compact fields
such as finite extensions of the $p$-adic rationals $\Qp.$
We also wish to lay some of the groundwork for a wavelet theory
in the adelic setting, with the goal of establishing another
technique from harmonic analysis for certain number theoretic
problems formulated in terms of adeles.
To these ends, this paper is devoted to the construction of a
wavelet theory on local fields.


\subsection{Results}\label{sbsec1.3}

Let $G$ be a LCAG
with compact open subgroup $H \subseteq G.$ 
Besides the cases that $G$ is compact or discrete,
there are many standard interesting examples of groups
$G$ with compact open subgroups; see Subsection~\ref{ssec:ex}
below, as well as \cite{BenR03}.
We develop a constructive theory of wavelets in $L^2(G)$
by means of the theory of wavelet sets mentioned in Subsection~\ref{sbsec1.1}.
In fact, we formulate a computable algorithm for generating
wavelet orthonormal bases (ONBs) for $L^2(G)$ depending on given
expansive automorphisms $A$ and given easily constructible mappings
$T$. We envisage that $A$ and $T$ will be chosen so that the resulting
wavelet ONB will have properties desirable for a given problem.

The algorithm is stated in Subsection~\ref{ssec:const}.
The proof that the algorithm does in fact generate wavelet
ONBs is given by Theorem~\ref{thm:wvltset} in Subsection~\ref{ssec:algpf}.
Theorem~\ref{thm:wvltset} depends on Theorem~\ref{thm:ONB} in
Subsection~\ref{ssec:thm1}. Theorem~\ref{thm:ONB} is geometrical in
nature; and it proves the equivalence of the existence of wavelets in terms
of sets $\Omega$ that are both ``translation'' and ``dilation'' tiles for
the dual group $\widehat{G}$.

The geometry of $G$ is decidedly non-Euclidean, and $G$ may not contain
a suitable discrete subgroup to serve as a lattice for translation. To
circumvent this problem we introduce the idea of $(\tau,\calD)$-congruence
in Subsections~\ref{ssec:trans} and \ref{ssec:wvlt}. Further, a notion
of expansive automorphism is required to define dilation properly.
This notion and its properties are formulated and proved in
Subsection~\ref{ssec:dil}. Subsection~\ref{ssec:ex} is devoted to
number theoretic examples and their behavior under our new ``translation''
and ``dilation'' operations.

Finally, Section~\ref{sec3} provides examples of wavelet ONBs for $L^2(G)$
generated by our algorithm. In particular, we prove 
in Proposition \ref{prop:haarshan} that the analogues of
the Euclidean Haar and Shannon wavelets are one and the same over $G$.
This section also gives us an opportunity to compare our method and its
generality with results of Lang \cite{Lan96}, \cite{Lan98a}, \cite{Lan98b},
which preceded this work, and Kozyrev \cite{Koz02}, which was brought to our
attention after our main results were proved, announced, and
written in a preliminary version of this paper. Because of the
more restrictive goals of Lang and Kozyrev, our results are significantly
more general than theirs, and our methods are fundamentally different.

Since we have to use the notion of an expansive automorphism, we 
hasten to reaffirm that our method is not classical MRA reformulated for $G$.
Our MRA theory of $L^2(G)$ is presented in \cite{BB04}. 

For all of its effectiveness in the Euclidean setting, MRA wavelet ONBs
require several steps of verification (e.g., using expansive linear
transformations to prove the density of $\bigcup V_j$ for an MRA
$\{V_j\}$) not necessary in our algorithm. On the other hand, the inherent
complexity of a wavelet set $\Omega$ and the slow wavelet decay of our
wavelets $\check{\mathbf 1}_{\Omega}$ on ${\mathbb R}^d$ are not a major
issue for $G$, since ${\mathbf 1}_{H^{\perp}}$ is locally constant,
that is, every point $\gamma\in \widehat{G}$ has
a neighborhood $U$ on which ${\mathbf 1}_{H^{\perp}}$ is constant.
This is true because $H^{\perp}$ is both an open and a closed
subset of $\widehat{G}$.  
($\check{\mathbf 1}_{\Omega}$ is the inverse Fourier transform of the
characteristic function 
${\mathbf 1}_{\Omega}$ of $\Omega$ and $H^{\perp}$ is the annihilator of
$H$; see Subsection~\ref{sbsec1.2}.)
Thus, ${\mathbf 1}_{H^{\perp}}$ is
locally smooth (meaning that ${\mathbf 1}_{H^{\perp}}$ would be
``$C^{\infty}$'' when there is a suitable
notion of derivative of a map from $\widehat{G}$ to $\CC$)
with compact support. Its inverse Fourier transform is
${\mathbf 1}_{H},$ which is also locally constant and hence is
also ``$C^{\infty}$'' with compact support.
Also, for a large class of the constructible mappings $T$, the
approximants to
${\mathbf 1}_{\Omega}$ at any finite step 
of our algorithm are 
locally constant with compact support, as are the
inverse transforms of such approximants.
This inverse Fourier transform generates a wavelet frame;
see~\cite{BS03} for the Euclidean case. Of course, the limit
${\mathbf 1}_{\Omega}$ itself is not necessarily locally constant,
although it does have compact support.

\subsection{Prerequisites and notation}\label{sbsec1.2}

Let $G$ be a LCAG, and denote its dual group by $\widehat{G}.$ 
Let $H$ be a closed subgroup of $G.$  
$$
     H^{\perp} = 
      \{\gamma \in \widehat{G}: \forall x \in H, \quad (x,\gamma) = 1\} 
     \subseteq \widehat{G}
$$
is the {\em annihilator subgroup} of $H$ in $\widehat{G}$,
where $(x,\gamma)$ denotes the
action of the duality between $G$ and $\widehat{G}.$

$\widehat{G/H}$ is algebraically and topologically isomorphic to 
$H^{\perp},$ and $\widehat{H}$ is algebraically and topologically isomorphic to
$\widehat{G}/H^{\perp}.$
$H^{\perp}$ is compact if and only if $H$ is open in $G$; and
$H$ is compact if and only if $H^{\perp}$ is open in $\widehat{G}$.
We shall suppress ``algebraically and topologically
isomorphic'' and write, for example, ``$\widehat{G/H}$ is $H^{\perp}.$'' 


$\mu=\mu_{G}$ and $\nu=\nu_{\widehat{G}}$
will denote Haar measures on $G$ and 
$\widehat{G},$ respectively. $L^2(G)$ is the space
of square integrable functions on $G,$
and the {\em Fourier transform} of $f \in L^2(G)$ is the function
$\hat{f} \in L^2(\widehat{G})$ formally defined as
$$
    \hat{f}(\gamma) = \int_G f(x)\bar{(x,\gamma)}\,d\mu(x),
$$
where $\gamma \in \widehat{G}$.

If the subgroup $H \subseteq G$ is open (and hence closed) and
compact, then $H$ and $H^{\perp}$ are compact abelian groups,
and the quotients $G/H$ and $\widehat{G}/H^{\perp}$ are
discrete abelian groups.  In that case, moreover,
we can and shall make
the following choices of normalization of Haar measure on
each of the six groups:
\begin{itemize}
\item $\mu$ satisfies $\mu(H)=1$,
\item $\nu$ satisfies $\nu(H^{\perp})=1$,
\item $\dsps \mu_{H} = \mu\bigr|_{H}$,
\item $\dsps \nu_{H^{\perp}} = \nu\bigr|_{H^{\perp}}$,
\item $\dsps \mu_{G/H}$ is counting measure, and
\item $\dsps \nu_{\widehat{G}/H^{\perp}}$ is counting measure.
\end{itemize}
These choices 
guarantee that the Fourier transform is an isometry
between $L^2(G)$ and $L^2(\widehat{G})$, and similarly between
$L^2(H)$ and $L^2(\widehat{H})=L^2(\widehat{G}/H^{\perp})$,
and between $L^2(G/H)$ and $L^2(\widehat{G/H})=L^2(H^{\perp})$,
e.g., \cite[\S 31.1]{HR70}, \cite{Rei68}. The above normalizations
are also compatible in the sense of Weil's formula \eqref{eq:weil},
see Remark \ref{rem:weil}.

Extra attention about normalization is necessary in the verification of the
aforementioned isometries in the case that $G$ is compact or discrete,
see \cite[\S 31.1]{HR70} again; but this case does not arise herein,
as noted in the discussion following Definition~\ref{def:expan}.
On the other hand, we do have to deal with the case $G = \widehat{G}$;
and some care is required, because the group
$\widehat{G/H}$ may not coincide
with the image of $H$ under the isomorphism $G\rightarrow \widehat{G}$.
For example, if $G$ is the self-dual
additive group of a finite extension $K$ of
$\Qp,$ then $H^{\perp} = \widehat{G/H}$ may not be $H=\ints_K$, but rather the
inverse different of $K;$ see, for example, \cite[\S 7.1]{RV99}.

We shall sometimes suppress $\mu$ and $\nu$ in lengthy
calculations; for example, we shall write ``$d\gamma$'' instead of
``$d\nu (\gamma)$'' in integrals. 

Let $\Aut G$ denote the group under composition of homeomorphic automorphisms
of $G$ onto itself.  For any $A\in\Aut G$ and for any Borel set
$U\subseteq G,$ $\mu_A(U)=\mu(AU)$
defines a Haar measure on $G.$ Formally,
$$
    \int_G f\circ{A^{-1}}(x)\,d\mu(x) = \int_G f(x)\,d\mu_A(x).
$$
Furthermore, $\mu_A$ cannot be the zero measure
because $\mu(U)=\mu_A(A^{-1}U)$.  Thus, there is a
unique positive number $|A|,$ the so-called {\em modulus} of
$A,$ with the property that
for any measurable set $U\subseteq G$, we have $\mu_A(U)=|A|\mu(U)$.
We shall use the following 
well-known fact, see \cite[\S 15.26]{HR63}:
\begin{equation}
\label{eq:CofV}
  \forall A \in \Aut G \; {\rm and} \; \forall f \in C_c(G), \quad
  \int_G f\circ A(x)\,dx = |A|^{-1}\,\int_G f(x)\,dx.
\end{equation}
Also, every 
element $A\in\Aut G$ has an adjoint element
$A^{\ast}\in\Aut \widehat{G}$, defined by
$(Ax,\gamma) = (x,A^{\ast}\gamma )$ for all $x\in G$ and
$\gamma\in\widehat{G}$.
We have
\begin{itemize}
\item $(A^{\ast})^{-1} = (A^{-1})^{\ast}$,
\item $|A|^{-1}=|A^{-1}|$, and 
\item $|A^{\ast}|=|A|$,
\end{itemize}
where the first equality is clear, the second
is immediate from \eqref{eq:CofV}, and the
third is by \eqref{eq:CofV} again and the Plancherel theorem.
Finally, and easily, the modulus function is a
homomorphism $\Aut G \longrightarrow (0,\infty)$,
where $(0,\infty)$ is considered as a multiplicative group.

Let $\calA\subseteq \Aut G$ be a nonempty countable
subset of $\Aut G$. In Theorem~\ref{thm:wvltset},
$\calA$ will be 
isomorphic to ${\mathbb Z}$;  however, we can phrase
the definitions and some results without assuming such a group
structure.

For any $x\in G$, we write $[x]$ or $x+H$ for its image in $G/H$;
that is, $[x]$ is the coset
\begin{equation}
\label{eq:coset}
[x] = \{ y\in G : y-x\in H\} \in G/H.
\end{equation} 
Similarly, for $\gamma$ in $\widehat{G}$, we write $[\gamma]$ 
or $\gamma + H^{\perp}$ for its image
in the quotient group $\widehat{G}/H^{\perp}.$ If $([x],\gamma )_{G/H}$
denotes the action of the duality between $G/H$ and $H^{\perp}$, then
$(x,\gamma )_G = ([x],\gamma )_{G/H}$ for $x \in G$, $\gamma\in H^{\perp}$.


These and other standard
facts about LCAGs that we shall use may be found in
\cite{HR63},
\cite{HR70},
\cite{Pon66}, 
\cite{Rei68},
and \cite{Rud62}. The standard
facts from number theory that we shall use are found in
\cite{RV99} and \cite{Rob00},
as well as \cite{Gou97}, \cite{Kob84}, and \cite{Ser79}.
\cite{Dau92} and \cite{Mey90}
are treatises on wavelet theory.

\begin{remark}
\label{rem:weil}
Classical Euclidean uniform sampling formulas depend essentially
on periodization in terms of a discrete subgroup.
Periodization induces a transformation from the given group to a
compact quotient group,
thereby allowing the analysis to be conducted in terms
of Fourier series which lead to sampling formulas.
The Shannon wavelet, which is a topic of Section~\ref{ssec:Haar},
is associated with the simplest (and slowly converging) Classical 
Sampling Formula derived from the $sinc$ function. {\em Weil's formula},
\begin{equation}
\label{eq:weil}
     \int_G f(x)\,d\mu(x) =
     \int_{G/H}\left( \int_H f(x+y)\,d\mu_H(y)\right)\,d\mu_{G/H}(x), 
\end{equation}
is a far reaching generalization of the idea of periodization,
which itself is manifested
in the term $\int_H f(x+y)\,d\mu_H(y).$ If two of the three Haar measures in
\eqref{eq:weil} are given then the third can be normalized so that
\eqref{eq:weil} is true on 
the space $C_c(G)$ of continuous functions with compact support.  
In our setting, with
$\mu_H$ as the restriction of $\mu$ to $H$, the choice
of $\mu_{G/H}$ to be counting measure is the appropriate
normalization for \eqref{eq:weil}. The analogous statement also applies
to $\nu$, $\nu_{H^{\perp}}$, and $\nu_{\widehat{G}/H^{\perp}}$.
\end{remark}

\section{Wavelets for groups with compact open subgroups}\label{sec:main}

\renewcommand{\theequation}{\thesection.\arabic{equation}}
\setcounter{equation}{0}

\subsection{Translations}
\label{ssec:trans}
Let $G$ be a LCAG
with compact open subgroup $H \subseteq G.$ 
$G$ may or may not contain a suitable discrete subgroup to
serve as a lattice for translations, e.g., Example~\ref{ex:Qp}.
One possible replacement for such a lattice is to choose a discrete set
$\calC$ of coset representatives for the quotient $G/H$,
and then translate elements of $L^2(G)$ by elements of $\calC$.
Kozyrev \cite{Koz02} used this strategy to 
show, for a certain choice $\calC$ of coset representatives
for $\Qp/\Zp$, that certain analogues of Haar wavelets
are wavelets in $L^2(\Qp)$ with respect to translations
by elements of $\calC$.

However, Kozyrev's examples
are for the group $\Qp$ only,
see the discussion after Proposition \ref{prop:haarshan}.  Moreover, 
there are substantial
obstacles to generalizing Kozyrev's method
to produce other wavelets, even for $\Qp$.  The main problem is
that because $\calC$ is not a group, there is no object
to serve as the dual lattice in $\widehat{\QQ}_p$.  Thus,
we shall need new translation operators on $L^2(G)$
which actually form a group.

Rather than choose elements $\calC$ in $G$,
we wish to construct one operator
$\tau_{[s]}$ for each element $[s]$ of
the discrete quotient group $G/H$, in such a way
that $\tau_{[s]+[t]}=\tau_{[s]}\tau_{[t]}$, and so that
$\tau_{[s]}$ is somehow similar to translation by $s$.
Our operators will be determined by a choice $\calD$ of 
coset representatives in $\widehat{G}$ for $\widehat{G}/H^{\perp}$.
(See Remark~\ref{rem:coset1} for more on $\calD$.)

\begin{definition}
\label{def:DetaW}
Let $G$ be a LCAG with 
compact open subgroup $H\subseteq G.$
Let $\calD\subseteq\widehat{G}$ be a set of coset
representatives in $\widehat{G}$ for the quotient
$\widehat{H}=\widehat{G}/H^{\perp}$.
\begin{list}{\rm \alph{bean}.}{\usecounter{bean}}
\item
Define maps
$\theta=\theta_{\calD}: \widehat{G} \rightarrow \calD$
and
$\eta=\eta_{\calD}: \widehat{G} \rightarrow H^{\perp} \subseteq \widehat{G}$
by
\begin{align}
\theta(\gamma) &= \mbox{ the unique } {\sigma}_{\gamma} \in\calD
	\mbox{ such that }
	\gamma-\sigma_{\gamma} \in H^{\perp} \subseteq \widehat{G},
\notag \\
\eta(\gamma) &= \gamma - \theta(\gamma).
\label{eq:etadef}
\end{align}

\item
For any fixed $[s]\in G/H$, define the following
unimodular (hence $L^{\infty}$) {\em weight function}:
\begin{equation}
\label{eq:wdef}
        w_{[s]}(\gamma) = w_{[s],\calD}(\gamma) =
	\overline{(s , \eta_{\calD}(\gamma))}.
\end{equation}
(See Remark~\ref{rem:wdef} for more on $w_{[s],\calD}$.)

\item
For any fixed $[s]\in G/H$, define the multiplier $m_{[s]}$
on $L^2(\widehat{G})$ as multiplication by $w_{[s],\calD}$,
that is, for any $F\in L^2(\widehat{G})$,
$$
        m_{[s]}F(\gamma) =
        m_{[s],\calD}F(\gamma) =
	F(\gamma) w_{[s],\calD}(\gamma) .
$$

\item
For any fixed $[s]\in G/H$ and for any $f\in L^2(G)$,
define $\tau_{[s]}f$ to be the inverse transform of
$m_{[s]}\widehat{f}$, that is,
$$
        {\tau}_{[s]}f =
        {\tau}_{[s],\calD}f =
	f \ast \check{w}_{[s],\calD},
$$
where the pseudo-measure $\check{w}_{[s],\calD}$ is
the inverse Fourier transform of $w_{[s],\calD}$.

\end{list}

\end{definition}

\begin{remark}
\label{rem:coset1}
a.
Recall that a set of coset representatives in $\widehat{G}$
for the quotient $\widehat{H}=\widehat{G}/H^{\perp}$ means
a set $\calD \subseteq \widehat{G}$
consisting of one element $\sigma \in \widehat{G}$
from each coset $\Sigma \in \widehat{G}/H^{\perp},$ i.e.,
$\Sigma = \sigma + H^{\perp}$.  For example,
if $\widehat{G}=\RR$ and $H^{\perp}=\ZZ$,
then the points of the
interval $[0,1)\subseteq\RR$ are a choice of coset representatives
for the quotient $\RR/\ZZ$.

b.
In our setting, the quotient
$\widehat{H}=\widehat{G}/H^{\perp}$ is discrete, and therefore
$\calD$ is also a discrete subset of $\widehat{G}$.  However,
just as $[0,1)$ is not a subgroup of $\RR$, our $\calD$ is not
necessarily a subgroup of $\widehat{G}$.  In fact, groups such
as $\widehat{G}=\Qp$ (with compact open subgroup $H^{\perp}=\ZZ_p$)
have no nontrivial discrete subgroups, and therefore cannot
have a $\calD$ which forms a group.

On the other hand, the image $[\calD] = \{[\sigma]: \sigma \in \calD\}$
of $\calD$ in the quotient $\widehat{G}/H^{\perp}$ is the full
group $\widehat{G}/H^{\perp} = \widehat{H}$, because $\calD$ contains
one element from {\em every} coset.

c.
Just as the choice of $[0,1)$ as a subset of $\RR$ is no
more natural than the choice of $(-1/2, 1/2]$, there 
is usually not a natural or canonical choice of the set $\calD$.
Nevertheless, in most cases, we can construct a set $\calD$ without
using the axiom of choice, e.g., Subsection~\ref{ssec:ex}
and Lemma~\ref{lem:Ddef} below.
\end{remark}

\begin{remark}
\label{rem:wdef}
a.
As the notation suggests, $w_{[s],\calD}$ depends
on $\calD$ and on the coset $[s]\in G/H$, but it does not depend on the
particular element $s\in G$.  Indeed, if $s,t\in G$, then
$$w_{[s+t],\calD}(\gamma) = \overline{(s + t , \eta(\gamma))}
= \overline{(s, \eta(\gamma))}\overline{(t , \eta(\gamma))}
= w_{[s],\calD}(\gamma) w_{[t],\calD}(\gamma).$$
The claim that $w_{[s],{\calD}}$ is independent of the particular
choice of $s$ in a given coset now follows from the observation
that if $t\in H$, then $w_{[t],{\calD}} = 1$.
It is also immediate that our operators have the desired
property that $\tau_{[s]+[t]}=\tau_{[s]}\tau_{[t]}$.

b.
The inverse Fourier transform of $w_{[s],\calD}$
is a pseudo-measure $\check{w}_{[s],\calD}$
on $G$; see \cite{Ben75} for pseudo-measures which are defined as
tempered distributions whose Fourier transforms are elements of
$L^{\infty}(\widehat{G})$.
In particular, the convolution $f \ast \check{w}_{[s],\calD}$
in the definition of $\tau_{[s],\calD}$
is well-defined by the action of the
Fourier transform, viz., 
$$
    \widehat{f \ast \check{w}}_{[s],\calD}
    = \widehat{f} w_{[s],\calD},
$$
which, in turn, is well-defined since 
$\widehat{f} w_{[s],\calD} \in L^2(\widehat{G}).$
Note that the operators $m_{[s],\calD}$
and $\tau_{[s],\calD}$ are unitary because $w_{[s],\calD}$
is unimodular.
\end{remark}

\begin{remark}
In our analogy with wavelets on $\RR^d$, $[s]$ takes the
place of a lattice element $n.$ Similarly, 
the pseudo-measure $\check{w}_{[s],\calD}$
takes the place of the $\delta$ measure at $n$.
Indeed, the $\delta$ measure at $s\in G$ acts
on $L^2(\widehat{G})$ by multiplication by
$\bar{(s,\gamma)}$, while $m_{[s]}$
acts by multiplication by 
$\bar{(s,\eta(\gamma))}$.  The ratio of the
two is $(s , \theta(\gamma) )$,
which has absolute value $1$ everywhere.

In the Euclidean setting, $G/H$ and $\calD$ are
both analogous to the discrete subgroup $\ZZ^d$, so that
both $s$ and $\theta(\gamma)$ would be elements
of $\ZZ^d$.  In particular, the $(s,\theta)$
factor would be identically $1$.
Thus, we may consider our operators $\tau_{[s]}$
to be as legitimate an analogue of the classical
Euclidean translation operators as are the
naive translation-by-$s$ operators.

In fact, if $s \in H,$ then ${w}_{[s],\calD} = 1,$ so that 
$\check{w}_{[s],\calD}$ is the $\delta$ measure at $0$,
and therefore
$f \ast \check{w}_{[s],\calD} = f$.  Generally, for
$s \notin H,$ $\check{w}_{[s],\calD}$ is more complicated,
but the action of $\tau_{[s]}$ is still easily computable,
see \cite{BenR03}.
\end{remark}

\subsection{Expansive automorphisms and dilations}
\label{ssec:dil}

When constructing wavelets in $L^2(\RR^d)$, one would like
an automorphism
(such as multiplication-by-$2$) which is expansive.
This means that the
automorphism must not only map the lattice into itself, but it
must have all eigenvalues greater than one in absolute value.
Although our group $G$ may have neither a lattice nor
a suitable notion of eigenvalues of automorphisms,
we can state a reasonable analogue of expansiveness as follows.

\begin{definition}
\label{def:expan}
Let $G$ be a LCAG with compact open
subgroup $H\subseteq G$, and let $A\in\Aut G$.  We say
that $A$ is {\em expansive} with respect to $H$ if both
of the following conditions hold:
\begin{list}{\rm \alph{bean}.}{\usecounter{bean}}
\item $H\subsetneq AH $, and
\item $\bigcap_{n\leq 0} A^n H = \{0\}$.
\end{list}
\end{definition}

The inclusion condition of expansiveness is the analogue of mapping
the lattice into itself as
in the definition of an MRA.
Similarly, the intersection condition is the analogue
of having eigenvalues greater than $1$, and a similar
condition in $\RR^d$ is a consequence of the other properties
in the definition of an MRA, e.g.,~\cite{Mad92}.

Note that if $A$ is an expansive automorphism of $G$,
then $|A|$ is an integer greater than $1$, just as is true of an
expansive integer matrix.
Indeed, $H$ is a proper subgroup of $AH$, so that $AH$ may be
covered by (disjoint) cosets $s+H$.  Each coset has measure $1$,
so that $AH$ must have measure equal to the number of cosets.
Thus, $|A|=\mu(AH)$ is just the
number of elements in the finite quotient group $(AH)/H$.

Further note that if $A$ is expansive, then $\mu(G)=\infty$.
To see this, simply observe
that by applying the first condition of Definition~\ref{def:expan}
repeatedly, we have
$$\mu(G)\geq \mu(A^n H) = |A|^n\mu(H) = |A|^n \geq 2^n$$
for every integer $n$.
As a result, $G$ cannot be compact, and $G/H$ is infinite.

We may alternatively characterize expansiveness by the action
of the adjoint automorphism $A^{\ast}$ on $\widehat{G}$,
as the following lemma shows.  Note that as a consequence
of this lemma, if $G$ has an expansive automorphism $A$,
then $\widehat{G}$ is $\sigma$-compact.  In addition,
$\widehat{G}$ cannot be compact, so $G$ cannot be discrete.
\begin{lemma}
\label{lem:aut}
Let $G$ be a LCAG with
compact open subgroup $H \subseteq G,$ and let $A \in \Aut G.$

a. $H \subseteq AH$ if and only if $H^{\perp} \subseteq A^{\ast}H^{\perp}.$

b. $H \subsetneq AH$ if and only if $H^{\perp} \subsetneq A^{\ast}H^{\perp}.$

c. Suppose $H \subseteq AH.$ Then
\begin{equation}
\label{eq:intuni}
     \bigcap_{n \leq 0}A^nH = \{0\} \Longleftrightarrow
     \bigcup_{n \geq 0}A^{{\ast}{n}}H^{\perp} = \widehat{G}.
\end{equation}
\end{lemma}

\begin{proof}
{\em a.}  Suppose $H\subseteq AH$, and consider $\gamma\in H^{\perp}$.
For any $x\in H$, we have $x\in AH$, and therefore $x=Ay$ for
some $y\in H$.  Hence,
$$(x, (A^*)^{-1} \gamma) = (Ay, (A^*)^{-1} \gamma)
=(y,\gamma)=1.$$
Thus, $(A^*)^{-1} \gamma \in H^{\perp}$, and so
$\gamma\in A^* H^{\perp}$.
It follows that $H^{\perp}\subseteq A^{\ast}H^{\perp}$.  The converse
is similar.

{\em b.}  If $H\subsetneq AH$, then $|A|>1$, and so
$A^*H^{\perp} \neq H^{\perp}$, because
$\nu(A^*H^{\perp})=|A^*|=|A|>1=\nu(H^{\perp})$.
Again, the converse is similar.

{\em c.} $(\Longleftarrow)$
Pick any
$x\in \bigcap_{n\leq 0} A^{n} H$, that is, any
$x\in G$ with $A^n x \in H$ for all $n\geq 0$.
For any $\gamma\in \widehat{G}$, there is an integer $n\geq 0$
such that $\gamma\in (A^*)^n H^{\perp}$, by assumption.
Thus,
$$(x,\gamma) = (A^n x , (A^*)^{-n} \gamma) = 1,$$
because $A^n x \in H$ and $(A^*)^{-n} \gamma\in H^{\perp}$.
Since $(x,\gamma)=1$ for every $\gamma\in \widehat{G}$, we must
have $x=0$.

{\em c.} $(\Longrightarrow)$
Pick any $\gamma\in \widehat{G}$, and suppose that
$\gamma\not\in \bigcup_{n\geq 0} (A^*)^n H^{\perp}$.
Then for every $n\geq 0$, there is some $x_n\in A^{-n} H$
such that $\zeta_n = (x_n,\gamma)\neq 1$.  Thus, there
is an integer $m\geq 1$ such that $\zeta_n^m$ has negative
real part.  Let $y_n = m x_n$, so that $y_n\in A^{-n} H$
(because $A^{-n} H$ is a group) and
$\text{Re}\, (y_n,\gamma) < 0$.

We claim that the sequence $\{y_n\}$ has an
accumulation point at $0 \in \widehat{G}$.  
Certainly some accumulation
point exists, because $y_n\in A^{-n} H \subseteq H$
and $H$ is compact.  If $0$ is not an accumulation point,
then there is a nonzero accumulation
point $y$.  Since $y\neq 0$, we have $y\not\in\bigcap A^{-n}H$, so
there must be an integer $N\geq 0$ for which
$y\not\in A^{-N}H$.  For any $n\geq N$, we have
$A^{-n} H \subseteq A^{-N}H$; therefore
$y\not\in A^{-n} H$.  By the previous paragraph,
then, $y + A^{-N} H$ is an
open neighborhood of $y$ containing only finitely
many $y_n$.  Hence, $y\neq 0$ is not actually an accumulation
point, and so the only accumulation point of $\{y_n\}$ is $y=0$.

Let $U=\{ x\in H : \text{Re}\, (x,\gamma) > 0\}$.
The map $(\cdot, \gamma) : G \rightarrow \CC$
is continuous by definition, and therefore $U$
is open.  Clearly $0\in U$.  Thus, there must be
infinitely many $y_n$ in $U$.  However,
$\text{Re}\, (y_n,\gamma) < 0$, and so $y_n\not\in U$.
That is a contradiction, and therefore $\gamma$
must lie in
$\bigcup_{n\geq 0} (A^*)^n H^{\perp}$.
\end{proof}

In Remark~\ref{rem:coset1}.c, it was noted that a set $\calD$
of coset representatives for an infinite quotient
$\widehat{G}/H^{\perp}$ can usually be constructed without using
the axiom of choice.  In fact,
given an expansive automorphism $A$ and a choice of coset
representatives for the finite quotient $H^{\perp}/((A^*)^{-1}H^{\perp})$,
we can generate such a $\calD$ by the following lemma.

\begin{lemma}
\label{lem:Ddef}
Let $G$ be a LCAG with compact
open subgroup $H\subseteq G$.
Let $A\in \Aut G$ be an expansive automorphism with modulus $|A|=N\geq 2$,
and let $\calD_1 = \{\rho_0,\ldots,\rho_{N-1}\}$ be
a set of coset representatives for the quotient
$H^{\perp}/((A^*)^{-1}H^{\perp})$, with $\rho_0=0$.
Define $\calD \subseteq \widehat{G}$ to be the
set of all elements $\sigma\in\widehat{G}$ of the form
\begin{equation}
\label{eq:Ddef}
\sigma = \sum_{j=1}^n (A^*)^j \rho_{i_j},
\quad \text{where } n\geq 1 \text{ and } i_j\in\{0,1,\ldots, N-1\}.
\end{equation}
Then $\calD$ is a set of coset representatives for the
quotient $\widehat{G} / H^{\perp}$.
Moreover, $A^* \calD\subsetneq \calD$.
\end{lemma}

\begin{proof}
First, we need to show that if $\sigma,\sigma'\in\calD$
with $\sigma-\sigma'\in H^{\perp}$, then
$\sigma=\sigma'$.  Write
$$\sigma = \sum_{j=1}^n (A^*)^j \rho_{i_j}, \quad
\sigma' = \sum_{j=1}^{n'} (A^*)^j \rho_{i'_j},$$
and assume without loss that $n\geq n'$.  Thus, we
have
\begin{equation}
\label{eq:sigsig}
\sum_{j=1}^{n'} (A^*)^j (\rho_{i_j}-\rho_{i'_j})
+ \sum_{j=n'+1}^n (A^*)^j \rho_{i_j} \in H^{\perp}.
\end{equation}
If $\sigma\neq\sigma'$, then
let $k$ be the largest value of $j$ for which the $(A^*)^j$
term is nonzero.  Note that every other term is
of the form either
$(A^*)^j (\rho_{i_j}-\rho_{i'_j}) \in
(A^*)^j H^{\perp}$
or $(A^*)^j \rho_{i_j} \in (A^*)^j H^{\perp}$,
with $j\leq k-1$.  Since
$(A^*)^j H^{\perp} \subseteq (A^*)^{k-1} H^{\perp}$,
we can multiply
both sides of the inclusion~\eqref{eq:sigsig}
by $(A^*)^{-k}$ to obtain
$$\rho_{i_k}-\rho_{i'_k} \in (A^*)^{-1} H^{\perp}$$
if $k\leq n'$, or
$$\rho_{i_k} \in (A^*)^{-1} H^{\perp}$$
if $k>n'$.  However, either of those conclusions implies that
the $(A^*)^k$ term of~\eqref{eq:sigsig} was zero, which contradicts
our assumption that $\sigma\neq\sigma'$.  Thus,
$\sigma=\sigma'$, as desired.

Second, given $\gamma\in\widehat{G}$, we need to show that
$\gamma\in\sigma+H^{\perp}$ for some $\sigma\in\calD$.
Let $n$ be the smallest nonnegative integer such that
$(A^*)^{-n}\gamma \in H^{\perp}$.  Such an integer must exist,
by Lemma~\ref{lem:aut}c.  We claim
that $\gamma$ is in $\sigma+H^{\perp}$ for some
$\sigma$ of the form
$$\sigma = \sum_{j=1}^n (A^*)^j \rho_{i_j}.$$
(Note that the $n$ in the sum is the same as the $n$ we
have just selected.)
We proceed by
induction on $n$.  If $n=0$, simply let $\sigma=0\in\calD$,
and we get $\gamma\in H^{\perp} = \sigma + H^{\perp}$.
For general $n\geq 1$,
let $\rho_{i_n}$ be the element of $\calD_1$ which
is in the same coset of $(A^*)^{-1}H^{\perp}$ as $(A^*)^{-n}\gamma$.
Then $(A^*)^{-(n-1)}(\gamma - (A^*)^n\rho_{i_n})\in H^{\perp}$,
and our induction is complete.

Finally, the fact that $A^*\calD\subseteq\calD$ is clear
from the form of $\sigma\in\calD$ and the fact that
$0\in\calD_1$; and $A^*\calD\neq\calD$ because
$A^{*}\rho_1\in\calD\setminus A^*\calD$.
\end{proof}

We define dilation operators $\delta_A$
(not to be confused with the $\delta$ measure)
in the usual way, as follows.
\begin{definition}
\label{def:dAdef}
Let $G$ be a LCAG, and let $A\in\Aut G$.
Define an operator $\delta_A$
on $L^2(G)$ by
\begin{equation}
\label{eq:dAdef}
\forall f\in L^2(G), \quad 
\delta_A f (x) = |A|^{1/2} f(Ax).
\end{equation}
\end{definition}

\begin{definition}
\label{def:diltrans}
Let $G$ be a LCAG with 
compact open subgroup $H\subseteq G$,
let $\calD$ be a choice of coset representatives in $\widehat{G}$
for $\widehat{H}=\widehat{G}/H^{\perp}$,
let $A\in\Aut G$,
and consider $[s]\in G/H$.
The {\em dilated translate} of $f\in L^2(G)$ is defined as
\begin{equation}
\label{eq:Asdef2}
f_{A,[s]}(x) = \delta_A \tau_{[s],\calD} f (x) =
|A|^{1/2} \cdot (f \ast \check{w}_{[s],\calD}) (Ax).
\end{equation}
We compute
\begin{equation}
\label{eq:Asdef}
\widehat{f_{A,[s]}}(\gamma) = |A|^{-1/2}
	\widehat{f}((A^{\ast})^{-1} \gamma)
	\overline{(s , \eta((A^{\ast})^{-1} \gamma) )}.
\end{equation}
\end{definition}

\subsection{Examples}
\label{ssec:ex}

The conditions that a LCAG $G$ contain
a compact open subgroup $H$ and have an expansive
automorphism $A$ place strong constraints on the
topologies of $G$ and $\widehat{G}$.
Such groups $G$ are a subclass of all totally disconnected groups.

\begin{figure}[hb]
\scalebox{.8}{\includegraphics{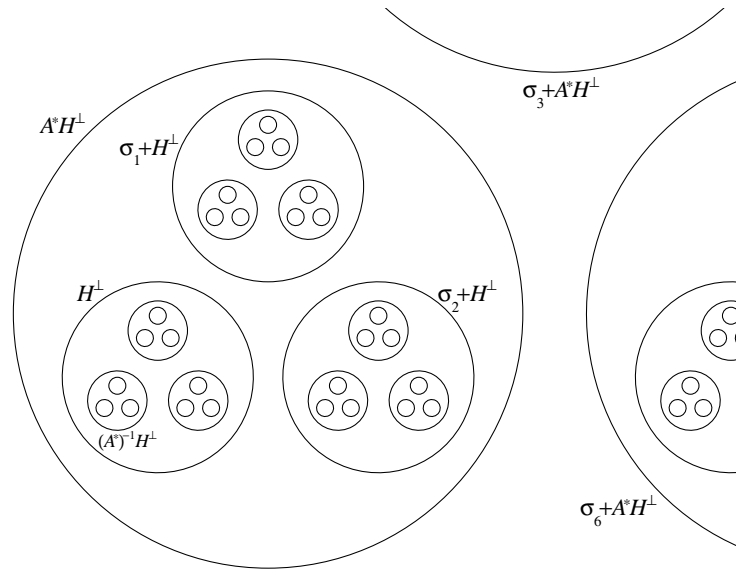}}
\caption{$\widehat{G}$ for a LCAG with compact open subgroup $H$
  and expansive automorphism $A$, with $|A|=3$.}
\label{fig:padic3}
\end{figure}

Figure~\ref{fig:padic3} indicates the structure of $\widehat{G}$
for these particular totally disconnected groups $G$ in the case
$|A|=3$. First note that $H^{\perp} \subsetneq A^{*}H^{\perp}$
by Lemma \ref{lem:aut}b. If ${\sigma}_1 \in A^{*}H^{\perp}\backslash H^{\perp}$
then $H^{\perp} \bigcap ({\sigma}_1 + H^{\perp}) = \emptyset$.
Continuing this procedure, we use a compactness argument to assert the
existence of finitely many cosets ${\sigma} + H^{\perp}$
covering $A^{*}H^{\perp}$, and this finite number is $|A|$ which is the
order of $A^{*}H^{\perp}/H^{\perp}$. Thus, we see that $A^{*}H^{\perp}$
is the finite union $H \bigcup ({\sigma}_1 + H^{\perp}) \bigcup 
({\sigma}_2 + H^{\perp})$ in Figure~\ref{fig:padic3}. The remaining
circles in Figure~\ref{fig:padic3} follow by considering
$A^{*}H^{\perp} \subsetneq {(A^{*})^2}H^{\perp}$ and
${(A^{*})^{-1}}H^{\perp} \subsetneq A^{*}H^{\perp}$. The idea of
the diagram is that the self-similar pattern of nested circles
continues forever.

$\widehat{G}$ is an infinite disjoint union of translates of 
$H^{\perp}$. In Figure~\ref{fig:padic3}, $H^{\perp}$ looks very
much like a Cantor set with constant rate of dissection determined
by $|A|$; thus, it is well to point out that $\widehat{G}$ need
not be totally disconnected.

In spite of the strong constraints mentioned above, many such groups $G$ 
and expansive automorphisms $A$ exist;
and many have applications in number theory.
We now give some examples
and comment on their translation and dilation operators.

\begin{example}
\label{ex:Qp}
Let $p\geq 2$ be a prime number.  The $p$-adic field $\Qp$
consists of all formal Laurent series in $p$ with coefficients
$0,1,\ldots,p-1$.  Thus,
$$
     \Qp = \left\{
     \sum_{n\geq n_0} a_n p^n : n_0\in\ZZ
     \text{ and } a_n \in \{0,1,\ldots, p-1\}
     \right\},
$$
with addition and multiplication as usual for Laurent series,
except with carrying of digits, so that, for example, in $\QQ_5$,
we have
$$(3 + 2\cdot 5) + (4 + 3 \cdot 5) = 2 + 1\cdot 5 + 1 \cdot 5^2.$$
(Equivalently, $\QQ_p$ is the completion of $\QQ$ with respect
to the $p$-adic absolute value $|p^r x|_p = p^{-r}$ for all $r\in\ZZ$
and all $x\in\QQ$ such that the numerator and denominator of $x$
are both relatively prime to $p$.)  Then $G=\Qp$ is a LCAG under
addition, with topology induced by $|\cdot|_p$, and with compact open
subgroup $H=\Zp$ consisting of Taylor series in $p$; equivalently,
$\Zp$ is the closure of $\ZZ\subseteq\Qp$.

$\Qp$ is self-dual, with duality action given by
$(x,\gamma)= \chi(x\gamma)$, where $\chi:\Qp \rightarrow \CC$
is the character given by
$$
\chi\left(\sum_{n\geq n_0} a_n p^n\right) =
\exp\left(2\pi i \sum_{n=n_0}^{-1} a_n p^n\right).
$$
The annihilator $\Zp^{\perp}$ is just $\Zp$ under this self-duality.

The quotient $\Qp/\Zp$ is isomorphic to $\mathbf{\mu}_{p^\infty}$,
the subgroup of $\CC^{\times}$ consisting of all roots of unity
$\zeta$ for which $\zeta^{p^n}=1$ for some $n\geq 0$.  One possible
choice for $\calD$, and the one used by Kozyrev \cite{Koz02}
in $G$ rather than in $\widehat{G}$, is
$$\calD_{\text{Koz}} = \left\{
\sum_{n=n_0}^{-1} a_n p^n : n_0\leq -1
\text{ and } a_n =0,1,\ldots, p-1
\right\},
$$
which is not a subgroup.
(Equivalently, $\calD_{\text{Koz}}$ is the set of all rational numbers
of the form $m/p^n$ for which $0\leq m\leq p^n -1$.)

A natural choice for $A$ would be multiplication-by-$1/p$, which
takes $\Zp$ to the {\em larger} subgroup $p^{-1}\Zp$; it
is easy to verify that this map is expansive.
More generally, for any nonzero $a\in\Qp$, the multiplication-by-$a$
map is an automorphism of the additive group $\Qp$,
with adjoint map also given by multiplication-by-$a$ on
the dual $\Qp$.  In fact, the only automorphisms of $\Qp$
are the multiplication-by-$a$ maps;
and these maps are expansive
if and only if $|a|_p>1$.
\end{example}

\begin{example}
\label{ex:Fpt}
Let $p\geq 2$ be a prime number, and let $\FF_p$ denote the
field of order $p$, with elements $\{0,1,\ldots, p-1\}$.
The field $\FF_p((t))$ consists of formal Laurent series in
the formal variable $t$, that is,
$$\FF_p((t)) = \left\{
\sum_{n\geq n_0} a_n t^n : n_0\in\ZZ
\text{ and } a_n \in \FF_p
\right\},$$
with the usual addition and multiplication of Laurent series,
but this time without carrying.
For example, in $\FF_5((t))$, we have
$$(3 + 2\cdot t) + (4 + 3 \cdot t) = 2.$$
(Equivalently, $\FF_p((t))$ is the completion of the field
$\FF_p(t)$ of rational functions with coefficients in $\FF_p$
with respect to the absolute value $|f|_0 = p^{-\ord_0(f)}$,
where $\ord_0(f)$ is the order of the zero (or negative the
order of the pole) of $f$ at $0$.)
Then $G=\FF_p((t))$ is a LCAG under addition, with the
topology induced by $|\cdot|_0$, and with compact open
subgroup $H=\FF_p[[t]]$ consisting of Taylor series in $t$;
equivalently, $H$ is the closure of the ring of polynomials
$\FF_p[t]\subseteq\FF_p((t))$.

$G$ is self-dual, with duality action given by
$(x,\gamma)= \chi(x\gamma)$, where $\chi:G \rightarrow \CC$
is the character given by
$$\chi\left(\sum_{n\geq n_0} a_n t^n\right) =
\exp\left(2\pi i a_{-1} p^{-1}\right).$$
The annihilator $H^{\perp}$ is just $H$ under this self-duality.

The quotient $G/H$ is isomorphic to a countable direct sum
of copies of $\FF_p$ (one copy for each negative power of $t$).
In fact,
$$\Lambda_{\FF_p((t))} = \left\{
\sum_{n= n_0}^{-1} a_n t^n : n_0\leq -1
\text{ and } a_n \in \FF_p
\right\}$$
is a discrete subgroup of $G$, and $G= H \times \Lambda_{\FF_p((t))}$
in the obvious way.  (The group $\Lambda_{\FF_2((t))}$ is the lattice
used by Lang \cite{Lan96,Lan98a,Lan98b} to study wavelets on $G=\FF_2((t))$,
which is sometimes called the Cantor dyadic group.)
$\Lambda_{\FF_p((t))}$, viewed as a subset of $\widehat{G}$,
is a logical choice for our $\calD$; in that case, because $\calD$
is actually a group, the operators $\tau_{[s],\calD}$ are in
fact the usual translation operators by elements of
$\Lambda_{\FF_p((t))}\subseteq G$, which is the annihilator
in $G$ of $\calD\subseteq\widehat{G}$.

One natural choice for $A$ would be multiplication-by-$1/t$, which
takes $H$ to the larger subgroup $t^{-1}H$.
As before, for any nonzero $a\in\Qp$, the multiplication-by-$a$
map is an automorphism of the additive group $G$, with adjoint
also multiplication-by-$a$; and it is expansive if and only
if $|a|_0>1$.  However, in contrast to $\Qp$,
many other automorphisms of $\FF_p((t))$ exist, because the coefficients
of the different powers of $t$ are arithmetically independent of
one another.  For example, one could define an automorphism
which acts (non-expansively)
only by multiplying the $t^{-17}$ coefficient by some
integer $m$ not divisible by $p$.
\end{example}

\begin{example}
\label{ex:L}
Let $L$ be any field which is a finite extension of either
$K=\Qp$ or $K=\FF_p((t))$.  Then the absolute value on $K$
extends in a unique way to $L$.  The ring of integers $\ints_L$
consists of all elements of $L$ which are roots of monic
polynomials with coefficients in $\ints_K$, where $\ints_K$
is $\Zp$ if $K=\Qp$ or $\FF_p[[t]]$ if $K=\FF_p((t))$.
Then $G=L$ is a LCAG under addition, with compact open subgroup
$H=\ints_L$.
We may choose an element $\pi_L\in \ints_L$ of largest possible
absolute value strictly less than $1$; such an element, called
a {\em uniformizer of $L$}, always exists, and it is analogous
to $t\in\FF_p((t))$ or $p\in\Qp$.  The quotient
$\ints_L/(\pi_L\ints_L)$ is a field of order $p^s$, for some
integer $s\geq 1$.

$G$ is self-dual, with duality action given by
$(x,\gamma)= \chi(\Tr(x\gamma))$, where $\chi:K \rightarrow \CC$
is as in the preceding examples, and $\Tr$ is the trace map from
$L$ to $K$.  This time, the annihilator $H^{\perp}$, known
as the {\em inverse different of $L$}, may be larger than $\ints_L$,
and it is of the form $\pi_L^{-r}\ints_L$ for some integer $r\geq 0$.

If $\calD_1$ is a list of $p^s$ coset representatives for the
finite quotient $\ints_L/(\pi_L\ints_L)$, then
$$\calD=\left\{
\sum_{n=n_0}^{-r-1} a_n \pi_L^n : n_0\leq -r-1
\text{ and } a_n \in \calD_1
\right\}
$$
is a set of coset representatives for $\widehat{G}/H^{\perp}$,
by Lemma \ref{lem:Ddef}.
If $K=\Qp$, then $\calD$ cannot be a group, as noted
in Example~\ref{ex:Qp};
but if $K=\FF_p((t))$, then we can choose $\calD_1$ to be a group,
in which case $\calD$ is also a group $\Lambda$,
as in Example~\ref{ex:Fpt}.

As before, for any nonzero $a\in L$, the multiplication-by-$a$
map is an automorphism of $L$ (expansive if and only if
$|a|>1$), with adjoint also given by
multiplication-by-$a$.  For extensions of $\Qp$, a few other
automorphisms involving Galois actions are possible; for
extensions of $\FF_p((t))$, many more automorphisms are
again possible.
\end{example}

\begin{example}
\label{ex:prod}
Let $G_1\ldots, G_N$ be LCAGs with compact open subgroups
$H_1,\ldots, H_N$ and coset representatives $\calD_i$ for
each quotient $\widehat{G}_i/ H^{\perp}_i$.  Then
$G=G_1\times\cdots \times G_N$ is a LCAG with compact open
subgroup
$H=H_1\times\cdots \times H_N$
and coset representatives
$\calD=\calD_1\times\cdots \times \calD_N$.
Thus, for example, $\Qp^N$ is an LCAG with compact open
subgroup, as is any product of fields from the preceding
examples.

An automorphism $A$ of $\Qp^N$ may be represented by
a matrix in $GL(N,\Qp)$, just as automorphisms of $\RR^N$
are matrices in $GL(N,\RR)$.  As is the case over $\RR$, the
adjoint automorphism is simply represented by
the transpose of the original matrix.  Such a matrix
map is expansive if and only if its eigenvalues have
absolute value ($|\cdot|_p$) greater than $1$.
\end{example}

\begin{example}
\label{ex:weird}
Let $X$ be a LCAG with compact open subroup $Y$, let
$W$ be a nontrivial discrete abelian group, and set $G = X \times W.$
Then $G$ is a LCAG with a compact open subgroup $H = Y \times \{0\}$.
Note that if $W$ is infinite, then,
even if $X$ is self-dual, the dual group
$\widehat{G}= \widehat{X} \times \widehat{W}$
need not be isomorphic to $G$.

Furthermore, the annihilator
$H^{\perp}$ is $Y^{\perp} \times\widehat{W}$.  Thus, if $A_X$
is an automorphism of $X$ which is expansive with respect to $Y$,
then $A=A_X\times \text{id}_W$ is an automorphism of $G$ which
is expansive with respect to $H$.  After all, the adjoint $A^{\ast}$
acts as $A_X^{\ast} \times \text{id}_{\widehat{W}}$, so that the
union of all positive iterates $(A^*)^n H^{\perp}$ is $\widehat{G}$.
This occurs in spite of the fact that the iterates $A^n H$ together
cover only $X \times \{0\}$, rather than all of $G$.
\end{example}

\begin{example}
\label{ex:aadic}
If $\mathbf{a}=\{a_n\}_{n\in\ZZ}\subseteq\NN$
is a bi-infinite sequence of positive
integers, one may define a group of $\mathbf{a}$-adics as formal
Laurent series in $t$ of the form
$$G=G_{\mathbf{a}} = \left\{
\sum_{n\geq n_0} c_n t^n : n_0\in\ZZ
\text{ and } c_n\in \{0,1,\ldots, a_n-1\}
\right\},$$
with usual addition of Laurent series, but
with carrying of digits, as in $\Qp$; see
\cite[\S 10.2]{HR63}.  (Thus, $\Qp$ is the
special case of $G_{\mathbf{a}}$ with $a_n=p$ for all $n$.)
Let $H=H_{\mathbf{a}}$ be the subgroup consisting
of Taylor series.  Then $G$ is a LCAG with compact open subgroup $H$.
As before, we could let $\calD$ consist of all Laurent series
for which all nonnegative power terms have coefficient $0$.

In this case, 
aside from uniform situations such as $\Qp$ for which
each $a_n$ is the same, $G$ generally
has no expansive isomorphisms, because the different powers of $t$
have incompatible coefficients.  Of course, one could
fix an integer $m$ and define $a_n=m$ for all $n$, so
that multiplication-by-$1/m$ would be an expansive automorphism
of the resulting group $G_{\mathbf{a}}$.
Then by the Chinese Remainder Theorem,
$G_{\mathbf{a}}$ would be
isomorphic to a product $\QQ_{p_1}\times\cdots\times\QQ_{p_s}$,
as $p_i$ ranges over the distinct prime factors of $m$;
such groups have therefore already been
discussed in the preceding examples.

\end{example}

\section{Theory of Wavelet Sets}

\subsection{Wavelets and $(\tau,\calD)$-congruence}
\label{ssec:wvlt}

Now that we have appropriate dilation and translation
operators, we are prepared to define wavelets on our
group $G$.

\begin{definition}
\label{def:newwvlt}
Let $G$ be a LCAG with compact open subgroup
$H \subseteq G$,
let $\calD\subseteq\widehat{G}$
be a choice of coset representatives in $\widehat{G}$
for $\widehat{H}=\widehat{G}/H^{\perp}$,
and let $\calA \subseteq \Aut G$ be a
countable nonempty set of automorphisms of $G$.
Consider $\Psi=\{\psi_1,\ldots,\psi_N\} \subseteq L^2(G).$ 
$\Psi$ is
a {\em set of wavelet generators} for $L^2(G)$ with respect to $\calD$
and $\calA$ if
$$
     \{\psi_{j,A,[s]} : 1\leq j \leq N, A\in\calA, [s]\in G/H \}
$$
forms an orthonormal basis for $L^2(G)$,
where
$$
\psi_{j,A,[s]}(x)=\delta_A\tau_{[s]}\psi_j(x)
= |A|^{1/2}\left(\tau_{[s],\calD} \psi_j \right) (Ax),$$
as in equation~\eqref{eq:Asdef2}.
In that case, the
resulting basis is called a {\em wavelet basis} for $L^2(G)$.

\indent
If $\Psi = \{\psi\},$
then $\psi$ is a {\em single wavelet} for $L^2(G).$
\end{definition}

The wavelets $\psi$ that we construct in this paper will have
the property that $\widehat{\psi}$ is the characteristic
function of some measurable subset of $\widehat{G}$.  We therefore
state the following definition, e.g.,~\cite{DL98,DLS97}.

\begin{definition}
\label{def:wvltset}
Let $G$, $H$, $\calD$, and $\calA$ be as in
Definition~\ref{def:newwvlt}.  Let $\Omega_1,\ldots,\Omega_N$ be
measurable subsets of $\widehat{G}$, and let
$\dsps \psi_j=\check{\mathbf{1}}_{\Omega_j}$ for each $j=1,\ldots, N$.
We say that $\{\Omega_1,\ldots,\Omega_N\}$ is a
{\em wavelet collection of sets}
if $\Psi=\{\psi_1,\ldots,\psi_N\}$ is a set of wavelet generators
for $L^2(G)$.

If $N=1$, then $\Omega=\Omega_1$ is a {\em wavelet set}.
\end{definition}

Based on a general formulation by Dai and Larson \cite{DL98} (with an
early version available in 1993) and by
Dai, Larson, and Speegle \cite{DLS97},
wavelet sets in $\widehat{\RR}^d$
were constructed 
by Soardi and Weiland \cite{SW98} and by Leon and one of the authors
\cite{BL99,BL01}
using the notion of $\tau$-congruence, which gives
an equivalence relation between subsets of $\widehat{\RR}^d$.
We now formulate an analogous equivalence relation for 
LCAGs $G$ with compact open subgroups $H.$

\begin{definition}
\label{def:semitau}
Let $G$ be a LCAG with compact open
subgroup $H \subseteq G$,
let $\calD\subseteq\widehat{G}$
be a choice of coset representatives in $\widehat{G}$
for $\widehat{H} = \widehat{G}/H^{\perp}$,
and let $V$ and $V'$ be subsets of $\widehat{G}$.
We say $V$ is
{\em $(\tau,\calD)$-congruent} to $V'$ 
if there is a finite or countable indexing set $I\subseteq\ZZ$
as well as
partitions $\{V_n : n\in I\}$ and $\{V'_n : n\in I\}$
of $V$ and $V'$, respectively, into measurable subsets,
and sequences
$\{\sigma_n\}_{n\in I}, \{\sigma'_n\}_{n\in I}\subseteq\calD$
such that
\begin{equation}
\label{eq:congru}
\forall n\in I, \qquad V_n\subseteq \sigma_n + H^{\perp}
\quad\text{and}\quad
V_n = V'_n -\sigma'_n + \sigma_n .
\end{equation}
\end{definition}

If $V_n = V'_n -\sigma'_n + \sigma_n$, then the condition that
$V_n\subseteq \sigma_n + H^{\perp}$ is equivalent
to $V'_n\subseteq \sigma'_n + H^{\perp}$.
Thus, in the case $V'=H^{\perp}$, we take
$\sigma'_n=\sigma'_0$ for all $n$,
where $\sigma'_0$ is the unique element
of $\calD\cap H^{\perp}$.
Note that the notion of $(\tau,\calD)$-congruence depends
crucially on the choice of coset representatives $\calD$.
Clearly, $(\tau,\calD)$-congruence is an equivalence relation,
and it preserves Haar measure.

\begin{remark}
a. Definition \ref{def:semitau} is related to bijective restrictions
of the canonical
surjection defined in \eqref{eq:coset}; see an analysis of this relation
in \cite[\S 3]{Ben98} in the context of Kluv{\'a}nek's sampling theorem (1965)
for LCAGs, cf., Remark~\ref{rem:weil}. Kluv{\'a}nek's
sampling formula
for a signal $f$ quantitatively relates the sampling rate with the measure of
the subsets of a given bandwidth corresponding to the frequency
content of $f.$

b. Congruence criteria were introduced
by Albert~Cohen (1990) to check the orthonormality of scaling functions 
of MRAs defined
by infinite products of dilations of a conjugate mirror filter, e.g.,
\cite{Dau92}, pages 182-186. The same notion of congruence also plays a
fundamental role in work on self-similar tilings by Gr{\"o}chenig, Haas,
Lagarias, Madych, Yang~B.~Wang, et al., e.g., \cite{GM92}, \cite{GH94},
\cite{LW95}, \cite{LW97}. We are dealing with (possibly infinite) sets
$\calD \subseteq \widehat{G}$ of coset representatives in $\widehat{G}$
for $\widehat{G}/H^{\perp}$, whereas the previous references deal with
coset representatives, so called {\em digits}, of each coset of $A({\Z}^d)$
in ${\Z}^d$. In this latter case, $A : {\R}^d \longrightarrow {\R}^d$
is an expansive matrix, and a goal is to determine those $A$ for which a
set of digits exists allowing for the construction of MRAs with 0-1
valued scaling functions.
\end{remark}

\subsection{Characterization of wavelet collections of sets}
\label{ssec:thm1}

Let $X$ be a measurable subset of $\widehat{G}$,
and let $\{X_n : n \in \ZZ \}$ be a countable set of measurable
subsets of $\widehat{G}$.  We say that
{\em $\{X_n\}$ tiles $X$ up to sets of measure zero}
if both
$$\nu\left( X \setminus
\left[ \bigcup_{n\in\ZZ} X_n \right] \right)= 0$$
and
$$
\forall m,n\in\ZZ, \, m\neq n, \quad
\nu(X_m \cap X_n) = 0.
$$
Similarly, we say that two sets $X,Y\subseteq\widehat{G}$ are
{\em $(\tau,\calD)$-congruent up to sets of measure zero}
if there are $(\tau,\calD)$-congruent sets
$X',Y'\subseteq\widehat{G}$
such that
$$\nu(X\setminus X') = \nu(X'\setminus X) =
\nu(Y\setminus Y') = \nu(Y'\setminus Y) = 0.$$
We are now prepared to state and prove our first main result,
which gives a necessary and sufficient condition
for a set $\Omega\subseteq\widehat{G}$ to be a wavelet set.

\begin{theorem}
\label{thm:ONB}
Let $G$ be a LCAG with compact open subgroup $H \subseteq G$,
let $\calD\subseteq\widehat{G}$
be a choice of coset representatives in $\widehat{G}$
for $\widehat{H}=\widehat{G}/ H^{\perp}$,
and let $\calA\subseteq\Aut G$ be a
countable nonempty subset of $\Aut G$.
Consider a finite sequence
$\{\Omega_1,\ldots,\Omega_N\}$ of measurable subsets of $\widehat{G}$.
$\{\Omega_1,\ldots,\Omega_N\}$ is a
wavelet collection of sets if and only if both
of the following conditions hold:
\begin{list}{\rm \alph{bean}.}{\usecounter{bean}}
\item
$\{A^{\ast}\Omega_j: A\in\calA, j=1,\ldots, N\}$ tiles $\widehat{G}$
up to sets of measure zero, and
\item
$\forall j=1,\ldots,N$,
$\Omega_j$ is $(\tau,\calD)$-congruent to $H^{\perp}$
up to sets of measure zero.
\end{list}
In that case, $\widehat{G}$ is $\sigma$-compact, each $\nu(\Omega_j)=1$,
and each $\mathbf{1}_{\Omega_j}\in L^2(\widehat{G})$.
\end{theorem}

\begin{proof}
Let $\psi_j=\check{\mathbf{1}}_{\Omega_j}$, and recall that
$\psi_{j,A,[s]} = \delta_A \tau_{[s],\calD} \psi_j$.
Recall also that $\sigma'_0$ is the unique 
element of $\calD\cap H^{\perp}$.

{\em i.}
We first verify that property {\em b}
implies that each $\nu(\Omega_j)=1$.
Let $I_j\subseteq\ZZ$ be the index set for the congruence from
Definition~\ref{def:semitau}, and let $\{V_{j,n} : n\in I_j\}$
be the corresponding partition of $\Omega_j$.
Since $\Omega_j$ is $(\tau,\calD)$-congruent to $H^{\perp}$,
we have
$$
\nu(\Omega_j)
= \sum_{n\in I_j} \nu(V_{j,n})
= \sum_{n\in I_j} \nu(V_{j,n} + \sigma'_n - \sigma_n)
= \sum_{n\in I_j} \nu(V'_{j,n})
= \nu(H^{\perp}) = 1,
$$
where $\{V'_{j,n} : n\in I_j\}$ is a partition of
$H^{\perp}$, and $\sigma'_n=\sigma'_0$ for all $n\in I_j$.
In particular, when property {\em b} is assumed,
we have $\mathbf{1}_{\Omega_j}\in L^2(\widehat{G})$.
Moreover, $\|\mathbf{1}_{\Omega_j}\|=1$,
so that $\|\psi_j\|_2 = 1$.

Conversely, if $\{\Omega_1,\ldots,\Omega_N\}$ is a wavelet collection
of sets, then
$$1=\|\psi_j\|_2 = \|\mathbf{1}_{\Omega_j}\|_2 = \nu(\Omega_j) .$$

We shall also need the fact that property {\em a}
implies that $\Omega_1,\ldots, \Omega_N$ are pairwise
disjoint up to sets of measure zero.  To see this,
pick any $A\in\calA$ and observe that for $j\neq k$,
$$
\nu(\Omega_j \cap \Omega_k)
=|A|^{-1} \nu\left( A^{\ast} (\Omega_j \cap \Omega_k) \right)
=|A|^{-1} \nu\left( (A^{\ast} \Omega_j) \cap (A^{\ast} \Omega_k) \right)
=0
$$
by property {\em a}.

{\em ii.}  We now prove that properties {\em a} and {\em b}
imply that
$\{\psi_{j,A,[s]} : 1\leq j\leq N, A\in\calA, [s]\in G/H \}$
is an orthonormal set.

First observe that
$\|\psi_{j,A,[s]}\|_2 = 1$
for any $A\in\calA$ and $[s]\in G/H$.
Indeed, we compute:
\begin{align*}
	\|\psi_{j,A,[s]}\|_2^2 &= 
	\int_{\widehat{G}} \left| \widehat{\psi_{j,A,[s]}}(\gamma) \right|^2
	d\, \gamma
	= |A|^{-1} \int_{\widehat{G}}
	\left| \widehat{\psi_j}((A^{\ast})^{-1} \gamma) \right|^2
	d\, \gamma  \\
	&= \int_{\widehat{G}}
	\left| \widehat{\psi_j}(\beta) \right|^2
	d\, \beta
	= \int_{\Omega_j}	d\, \beta
	= \nu(\Omega_j) = 1,
\end{align*}
by the change of variables $\beta=(A^{\ast})^{-1}\gamma$.
Furthermore, if $(j,B,[t])$ and $(k,C,[u])$ are distinct triples
in $\{1,\ldots, N\}\times \calA \times (G/H)$, then
$\langle\psi_{j,B,[t]} , \psi_{k,C,[u]}\rangle = 0.$ To see this, 
consider the cases $(j,B) \neq (k,C)$ and
$(j,B) = (k,C).$

If $(j,B) \neq (k,C),$ then we compute
\begin{align*}
& \langle\psi_{j,B,[t]} , \psi_{k,C,[u]}\rangle \\
& =
(|B| |C|)^{-1/2} \int_{\widehat{G}}
\widehat{\psi_j}\left( (B^{\ast})^{-1} \gamma\right)
\overline{\widehat{\psi_k}\left( (C^{\ast})^{-1} \gamma\right)}
\overline{\left( t,\eta\left( (B^{\ast})^{-1} \gamma\right) \right)}
\left( u,\eta\left( (C^{\ast})^{-1} \gamma\right) \right) d\gamma \\
& =
(|B| |C|)^{-1/2} \int_{\widehat{G}}
\mathbf{1}_{B^{\ast}\Omega_j}(\gamma)
\mathbf{1}_{C^{\ast}\Omega_k}(\gamma)
\overline{\left( t,\eta\left( (B^{\ast})^{-1} \gamma\right) \right)}
\left( u,\eta\left( (C^{\ast})^{-1} \gamma\right) \right) d\gamma
=0 ,
\end{align*}
because the intersection
$(B^{\ast}\Omega_j) \cap (C^{\ast}\Omega_k)$
has measure zero, by property {\em a}.

If $(j,B) = (k,C)$ but $[t]\neq [u]$, recall
that, by property {\em b},
$\Omega_j$ is $(\tau,\calD)$-congruent to $H^{\perp}$ up to sets
of measure zero.  We may therefore
write $\Omega_j$ (up to sets of measure zero)
as the disjoint union of sets $V_{j,n}$,
where $n$ ranges over a finite or countable index set $I_j$,
and where $H^{\perp}$ is the disjoint 
union of the sets $V'_{j,n}=V_{j,n} - \sigma_{j,n} +\sigma'_0$,
where each $\sigma_{j,n}\in\calD$.
Therefore,
\begin{align}
\label{eq:ut}
\langle\psi_{j,B,[t]} , \psi_{j,B,[u]}\rangle &=
\int_{\Omega_j} (u-t,\eta(\gamma)) d\gamma
\notag
\\
&= \sum_{n\in I_j} \int_{V_{j,n}}
(u-t,\eta(\gamma)) d\gamma
= \sum_{n\in I_j} \int_{V'_{j,n}-\sigma'_0 + \sigma_n}
(u-t,\eta(\gamma)) d\gamma
\notag
\\
&= \sum_{n\in I_j} \int_{V'_{j,n}}
(u-t,\eta(\lambda - \sigma'_0 + \sigma_n)) d\lambda
= \sum_{n\in I_j} \int_{V'_{j,n}}
(u-t,\lambda - \sigma'_0) d\lambda
\notag
\\
&= \overline{ (u-t,\sigma'_0) }
\int_{H^{\perp}} (u-t,\lambda) d\lambda,
\end{align}
where the second, third, and sixth equalities follow from the definition
of $(\tau,\calD)$-congruence, and the fourth is the change of variables
$\lambda = \gamma -\sigma_n + \sigma'_0$.
The fifth follows from the definition of $\eta$,
since $\lambda-\sigma'_0\in H^{\perp}$, so that
$\theta(\lambda - \sigma'_0 + \sigma_n) = \sigma_n$,
and hence
$$\eta(\lambda - \sigma'_0 + \sigma_n)
= (\lambda - \sigma'_0 + \sigma_n) - \sigma_{n}
= \lambda - \sigma'_0 .$$
Finally, by the following well known proof, \eqref{eq:ut} vanishes
because $[t]\neq [u]$, i.e., because $x=u-t\not\in H^{\perp}$:
$$\int_{H^{\perp}} (x,\lambda) d\nu_{H^{\perp}}(\lambda)=
(x,\gamma_0) \int_{H^{\perp}}
(x,\lambda - \gamma_0) d\nu_{H^{\perp}}(\lambda)=
(x,\gamma_0) \int_{H^{\perp}}
(x,\lambda) d\nu_{H^{\perp}}(\lambda),$$
where we have used the translation invariance of $\nu_{H^{\perp}}$
and where $\gamma_0\in H^{\perp}$ was chosen so that
$(x,\gamma_0)\neq 1$.

Thus, $\{\psi_{j,A,[s]} : 1\leq j\leq N, A\in\calA, [s]\in G/H \}$
is an orthonormal set.

{\em iii.}  Next, we prove that properties {\em a} and {\em b}
are sufficient for $\{\Omega_1, \ldots, \Omega_N\}$
to be a wavelet collection of sets.
The substantive part
of the proof is concerned with proving the following formula:
\begin{equation}
\label{eq:maincl}
        \forall f\in L^2(G), \quad
	\sum_{\substack{1\leq j\leq N \\ A\in\calA \\ [s]\in G/H} }
	\left| \langle f , \psi_{j,A,[s]} \rangle \right|^2
	= \| f \|_2^2.
\end{equation}
This fact, combined with the orthonormality, proves that
$\{\psi_{j,A,[s]} : 1\leq j\leq N, A\in\calA, [s]\in\ G/H \}$
is an orthonormal basis for $L^2(G).$

Now to the proof of \eqref{eq:maincl}. 
By the orthonormality, the left side of \eqref{eq:maincl} is 
bounded by $\| f \|_2^2$ when the sum is over any finite set of indices
(Bessel's inequality). Thus,
$$
     \forall f \in L^2(G), \quad {\rm card} \{\psi_{j,A,[s]}:
     \langle f , \psi_{j,A,[s]} \rangle
     \neq 0\} \leq \aleph_o, 
$$ 
and so the left side of \eqref{eq:maincl} is finite,
and we shall evaluate it.
By Plancherel's theorem
and \eqref{eq:Asdef}, we have
\begin{align}
	\sum_{\substack{1\leq j\leq N \\ A\in\calA \\ [s]\in G/H} }
	\left| \langle f , \psi_{j,A,[s]} \rangle \right|^2
	&=
	\sum_{\substack{1\leq j\leq N \\ A\in\calA \\ [s]\in G/H} }
	\left|
	\langle \hat{f} , \widehat{\psi_{j,A,[s]}} \rangle
	\right|^2
	\notag \\
	&=
	\sum_{\substack{1\leq j\leq N \\ A\in\calA \\ [s]\in G/H} }
	|A|^{-1} \left|
	\int_{\widehat{G}} \hat{f}(\gamma)
	\overline{\widehat{\psi_j}((A^{\ast})^{-1} \gamma)}
	(s , \eta((A^{\ast})^{-1} \gamma) )
	\, d\gamma \right|^2
	\notag \\
	&=
	\sum_{\substack{1\leq j\leq N \\ A\in\calA \\ [s]\in G/H} }
	|A| \left|
	\int_{\widehat{G}} \hat{f}(A^{\ast}\beta)
	\overline{\widehat{\psi_j}(\beta)}
	(s , \eta(\beta) )
	\, d\beta \right|^2 ,
	\label{eq:fprod1}
\end{align}
where we have substituted
$\beta= (A^{\ast})^{-1}\gamma,$ and hence
$d\beta = |A|^{-1} d\gamma,$ i.e., \eqref{eq:CofV}.

By property {\em b},
$\Omega_j$ is $(\tau,\calD)$-congruent to $H^{\perp}$ up to sets
of measure zero.
Recalling the notation $I_j$,
$V'_{j,n}$, $V_{j,n}$, and $\sigma_{j,n}$
from the discussion preceding \eqref{eq:ut},
the right side of \eqref{eq:fprod1} becomes
\begin{align}
        &
        \sum_{\substack{1\leq j\leq N \\ A\in\calA \\ [s]\in G/H} }
	|A| \left|
	\int_{\Omega_j} \hat{f}(A^{\ast}\beta)
	(s , \eta(\beta))
	\, d\beta \right|^2
	\notag \\
	= &
        \sum_{\substack{1\leq j\leq N \\ A\in\calA \\ [s]\in G/H} }
	|A| \left|
	\sum_{n\in I_j}
	\left[
	\int_{V_{j,n}} \hat{f}(A^{\ast}\beta)
	(s , \eta(\beta))
	\, d\beta \right] \right|^2
	\notag \\
	= &
        \sum_{\substack{1\leq j\leq N \\ A\in\calA \\ [s]\in G/H} }
	|A| \left|
	\sum_{n\in I_j}
	\left[
	\int_{V'_{j,n}-\sigma'_0} \hat{f}(A^{\ast}(\alpha + \sigma_{j,n}))
	(s , \eta(\alpha +\sigma_{j,n}))
	\, d\alpha \right] \right|^2 \notag \\
	= &
        \sum_{\substack{1\leq j\leq N \\ A\in\calA \\ [s]\in G/H} }
	|A| \left|
	\sum_{n\in I_j}
	\left[
	\int_{\widehat{G}}
	\mathbf{1}_{V'_{j,n}}(\alpha+\sigma'_0)
	\hat{f}(A^{\ast}(\alpha + \sigma_{j,n}))
	(s , \eta(\alpha +\sigma_{j,n}))
	\, d\alpha \right] \right|^2 ,
\label{eq:fprod2}
\end{align}
where we have substituted $\alpha = \beta -\sigma_{j,n}$.
Since $\alpha\in H^{\perp} - \sigma'_0 =H^{\perp},$ we
have
$$\eta(\alpha + \sigma_{j,n})
= (\alpha + \sigma_{j,n}) - \theta(\alpha +\sigma_{j,n})
= (\alpha + \sigma_{j,n}) - \sigma_{j,n} = \alpha.$$
Furthermore, we can
exchange the inner summation and integral signs
in the last term of \eqref{eq:fprod2}.  In fact,
since $\{V'_{j,n}\}$ is a tiling of $H^{\perp}$ and
$K_n = V'_{j,n} - \sigma'_0 \subseteq H^{\perp}$,
then $\{K_n\}$ is a tiling of $H^{\perp}$; and, hence,
denoting the integrand by $F_n$ (noting it vanishes
off $K_n$), and writing $F=\sum F_n$,
we see that $F_n, F \in L^2(H^{\perp}) \subseteq L^1(H^{\perp})$,
that $\bigcup K_n = H^{\perp}$, that
$$\sum_{n\in I_j} \int_{\widehat{G}} F_n(\alpha) d\alpha
= \sum_{n\in I_j} \int_{K_n} F(\alpha) d\alpha
= \int_{H^{\perp}} F(\alpha) d\alpha.$$
Thus, the right side of~\eqref{eq:fprod2} becomes
\begin{equation}
\label{eq:fprod3}
        \sum_{\substack{1\leq j\leq N \\ A\in\calA} }
	|A| \sum_{[s]\in G/H}
	\left|
	\int_{H^{\perp}}
	\left[
	\sum_{n\in I_j}
	\mathbf{1}_{V'_{j,n}}(\alpha+\sigma'_0)
	\hat{f}(A^{\ast}(\alpha + \sigma_{j,n}))
	\right]
	(s , \alpha )
	\, d\alpha \right|^2 .
\end{equation}

Let $F_{j,A}(\alpha) = \sum_{n\in I_j} \mathbf{1}_{V'_{j,n}}(\alpha+\sigma'_0)
\hat{f}(A^{\ast}(\alpha + \sigma_{j,n}))$ on $H^{\perp}\subseteq\widehat{G}.$ 
In particular, $\text{supp}\,F_{j,A} \subseteq H^{\perp}.$
We consider $F_{j,A}$ as an element of 
$L^2(H^{\perp}),$ where the Haar measure $\nu_{H^{\perp}}$
may be considered the restriction of $\nu = \nu_{\widehat{G}}$
to the compact group $H^{\perp} \subseteq \widehat{G}$
with the scaling $\nu_{H^{\perp}}(H^{\perp}) = \nu(H^{\perp}) = 1$,
as noted in Subsection~\ref{sbsec1.2}.
As such, we have
$$
\int_{\widehat{G}} F_{j,A}(\alpha)d\nu(\alpha)
=\int_{H^{\perp}} F_{j,A}(\alpha)d\nu_{H^{\perp}}(\alpha) .
$$
Hence, since $G/H$ is the 
discrete dual group of $H^{\perp}$, we may apply Plancherel's theorem
for Fourier series and obtain
$$\sum_{[s]\in G/H} \left|
	\int_{H^{\perp}} F_{j,A}(\alpha) ( [s], \alpha ) \, 
	d\nu(\alpha) \right| ^2 =
	\int_{H^{\perp}} \left| F_{j,A} (\alpha) \right|^2 \,
	d\nu(\alpha) .
$$
Thus, because $([s],\alpha)=(s,\alpha)$ for $[s]\in G/H$ and
$\alpha\in H^{\perp}$, \eqref{eq:fprod3} becomes
$$
	\sum_{\substack{1\leq j\leq N \\ A\in\calA} }
	|A| \int_{H^{\perp}}
	\left|
	\sum_{n\in I_j}
	\mathbf{1}_{V'_{j,n}}(\alpha+\sigma'_0)
	\hat{f}(A^{\ast}(\alpha + \sigma_{j,n}))
	\right|^2 \, d\alpha 
$$
which, in turn, is
\begin{equation}
\label{eq:fprod4}
	\sum_{\substack{1\leq j\leq N \\ A\in\calA} }
	|A| \int_{H^{\perp}}
	\left[
	\sum_{n\in I_j}
	\left|
	\hat{f}(A^{\ast}(\alpha + \sigma_{j,n}))
	\right|^2
	\mathbf{1}_{V'_{j,n}}(\alpha+\sigma'_0)
	\right]
	\, d\alpha ,
\end{equation}
where \eqref{eq:fprod4} is a consequence of the fact that,
for $j$ fixed, the sets $V'_{j,n}$ 
are pairwise disjoint, so that for any given $\alpha$
at most one term in the sum over $n$ is nonzero.
We can now interchange the inner summation and integral as before.
Hence,~\eqref{eq:fprod4} 
becomes
\begin{align}
        & \sum_{\substack{1\leq j\leq N \\ A\in\calA} }
        |A| \sum_{n\in I_j}
	\left[
	\int_{H^{\perp}}
	\left|
	\hat{f}(A^{\ast}(\alpha + \sigma_{j,n}))
	\right|^2
	\mathbf{1}_{V'_{j,n}}(\alpha+\sigma'_0)
	\, d\alpha  \right] \notag \\
	= & \sum_{\substack{1\leq j\leq N \\ A\in\calA} }
	|A| \sum_{n\in I_j}
	\left[
	\int_{V'_{j,n}-\sigma'_0}
	\left|
	\hat{f}(A^{\ast}(\alpha + \sigma_{j,n}))
	\right|^2
	\, d\alpha  \right] \notag \\
	= & \sum_{\substack{1\leq j\leq N \\ A\in\calA} }
	|A| \sum_{n\in I_j}
	\left[
	\int_{V_{j,n}}
	\left|
	\hat{f}(A^{\ast}(\beta))
	\right|^2
	\, d\beta  \right] \notag \\
	= & \sum_{\substack{1\leq j\leq N \\ A\in\calA} }
	|A| \int_{\Omega_j}
	\left|
	\hat{f}(A^{\ast}(\beta))
	\right|^2
	\, d\beta \notag \\
	= & \sum_{\substack{1\leq j\leq N \\ A\in\calA} }
	\int_{A^{\ast}\Omega_j}
	\left|
	\hat{f}(\gamma)
	\right|^2
	\, d\gamma .
\label{eq:fprod5}
\end{align}
However, $\{A^{\ast}\Omega_j\}$ tiles $\widehat{G}$
up to sets of measure zero. 
Hence, the right side of~\eqref{eq:fprod5} is
$$
\int_{\widehat{G}} \left| \hat{f}(\gamma) \right|^2 \, d\gamma
= \|\hat{f}\|_2^2 = \|f\|_2^2,
$$
and so formula~\eqref{eq:maincl} is proved.

{\em iv}.  Let $\{\Omega_1,\ldots,\Omega_N\}$ be a wavelet
collection of sets.  We now prove that $\widehat{G}$ is
$\sigma$-compact.  Note that for any
$A\in\calA$ and $1\leq j\leq N$,
$$\nu(A^{\ast}\Omega_j) = |A|\nu(\Omega_j) = |A| < \infty.$$
Thus, there are at most countably many
$\sigma\in\calD$ such that
$\nu((\sigma+H^{\perp}) \cap (A^{\ast}\Omega_j) ) >0 .$
Write 
$X=\bigcup_{A\in\calA, 1\leq j\leq N} A^{\ast} \Omega_j$,
which is a countable union; hence, there are at
most countably many $\sigma\in\calD$ such that
$\nu((\sigma+H^{\perp}) \cap X ) > 0.$

Suppose $\calD$ were uncountable; then there is
some $\sigma\in\calD$ such that $\nu((\sigma+H^{\perp}) \cap X ) = 0$.
Let $F=\mathbf{1}_{\sigma+H^{\perp}}\in L^2(\widehat{G})$.
Then for any $A\in\calA$, $1\leq j\leq N$,
and $[s]\in G/H$, we have
$\langle F, \widehat{\psi_{j,A,[s]}}\rangle = 0$,
because $\widehat{\psi_{j,A,[s]}}$ vanishes off of $X$.
Therefore,
$$1=\|F\|_2^2 = \sum_{j,A,[s]}
|\langle F, \widehat{\psi_{j,A,[s]}}\rangle|^2=0.$$
By this contradiction, it follows that $\calD$ is countable.
Thus, $\widehat{G}=\bigcup_{\sigma\in\calD}(\sigma+H^{\perp})$
is a countable union of compact sets, as desired.

{\em v.}
Finally we prove that properties {\em a} and {\em b} are necessary
conditions for
$\{\psi_j\}$ to be a set of wavelet generators.

First we show that $\{A^{\ast}\Omega_j\}$ tiles $\widehat{G}$
up to set of measure zero (property {\em a}).
If $(j,B),(k,C)\in\{1,\ldots, N\}\times \calA$ are distinct pairs,
then $\langle \psi_{j,B,[0]},\psi_{j,C,[0]} \rangle=0$
by orthonormality.  Thus,
$$\nu(B^{\ast}\Omega_j\cap C^{\ast}\Omega_k) =
\langle \mathbf{1}_{B^{\ast}\Omega_j}, \mathbf{1}_{C^{\ast}\Omega_k} \rangle
=
\langle \psi_{j,B,[0]}, \psi_{k,C,[0]} \rangle = 0.$$
If $\alpha\in 
Y=\left( \widehat{G} \setminus
\left[ \bigcup_{A\in\calA, 1\leq j\leq N} A^{\ast} \Omega_j \right] \right)$,
then let $Y_{\alpha} = Y\cap (\alpha + H^{\perp})$.
Clearly $\mathbf{1}_{Y_{\alpha}}\in L^2 (\widehat{G})$,
and the orthonormality
of $\{\psi_{j,A,[s]}\}$ ensures that
${\rm card} \{(j,A,[s]) :
\langle \mathbf{1}_{Y_{\alpha}}, \psi_{j,A,[s]} \rangle \neq 0 \}
\leq \aleph_0$.  Since $\{\psi_{j,A,[s]}\}$ is an
orthonormal basis for $L^2(G)$, we have
\begin{align*}
\nu(Y_{\alpha}) & =
\| \mathbf{1}_{Y_{\alpha}}\|_2^2 = 
\sum_{\substack{1\leq j\leq N \\ A\in\calA \\ [s]\in G/H} }
\left| \langle
\mathbf{1}_{Y_{\alpha}},\widehat{\psi_{j,A,[s]}}
\rangle \right|^2
\notag \\
&=
\sum_{\substack{1\leq j \leq N \\ A\in\calA \\ [s]\in G/H} }
|A|^{-1} \left|
\int_{\widehat{G}} \mathbf{1}_{Y_{\alpha}}(\gamma)
\overline{\widehat{\psi_j}((A^{\ast})^{-1} \gamma)}
(s , \eta((A^{\ast})^{-1} \gamma) )
\, d\gamma \right|^2
\notag \\
&=
\sum_{\substack{1\leq j\leq N \\ A\in\calA \\ [s]\in G/H} }
|A|^{-1} \left|
\int_{Y_{\alpha}\cap A^{\ast}\Omega_j}
(s , \eta((A^{\ast})^{-1} \gamma) )
\, d\gamma \right|^2 = 0 ,
\end{align*}
because $Y_{\alpha}\cap A^{\ast}\Omega_j=\emptyset$.
We have covered $Y$ by open sets $\alpha + H^{\perp}$ such
that each intersection $Y\cap (\alpha + H^{\perp})$ has measure zero;
thus, $\nu(Y)=0$, by the $\sigma$-compactness of $\widehat{G}$.

Second, we show that each $\Omega_j$ is $(\tau,\calD)$-congruent
to $H^{\perp}$ (property {\em b}).  Since $\nu(\Omega_j)=1 < \infty$,
there are at most countably many cosets $\sigma + H^{\perp}$
such that $\Omega_j\cap(\sigma+ H^{\perp})$ has positive measure.
(In fact, $\widehat{G}/H^{\perp}$ is countable because
$\widehat{G}$ is $\sigma$-compact.)
Let $\{\sigma_{j,n} : n\in I_j\}$ be the set of all such $\sigma\in\calD$,
where $I_j\subseteq\ZZ$ is some finite or countable indexing set.
For each $n\in I_j$,
let $V_{j,n} = \Omega_j\cap (\sigma_{j,n} + H^{\perp})$
and let $V'_{j,n} = V_{j,n} - \sigma_{j,n} + \sigma'_0$,
where $\sigma'_0$ is
the unique element of $\calD\cap H^{\perp}$.

Clearly, $\{V_{j,n}\}$ is a partition of $\Omega_j$,
by the definitions of $\sigma_{j,n}$ and $V_{j,n}$.

It remains to prove that
$\{V'_{j,n}\}$ is a partition of $H^{\perp}$.
The proof requires our hypothesis of orthonormality,
as follows.
Let $j\in\{1,\ldots, N\}$ and $A\in\calA$ be
fixed, and compute, for any $[s]\in G/H$, that
\begin{align}
\label{eq:3.AA}
\delta([s]) & = \langle \psi_{j,A,[0]} , \psi_{j,A,[s]} \rangle
\notag
\\
& = |A|^{-1} \int_{\widehat{G}}
\left| \widehat{\psi}_j \left( (A^{\ast})^{-1} \gamma \right) \right|^2
\overline{ \left(s, \eta \left( (A^{\ast})^{-1} \gamma \right) \right)}
d\gamma
\notag
\\
& = \int_{\Omega_j} \overline{ (s, \eta(\gamma) )} d\gamma
= \sum_{n\in I_j} \int_{V_{j,n}} \overline{ (s, \eta (\gamma) )} d\gamma
\notag
\\
& = \sum_{n\in I_j} \int_{V_{j,n}-\sigma_{j,n}}
\overline{ (s, \eta(\sigma_{j,n} + \beta) )} d\beta
\notag
\\
& = \sum_{n\in I_j} \int_{V_{j,n}-\sigma_{j,n}}
\overline{ (s, \beta )} d\beta ,
\end{align}
where $\delta([s])=0$ if $[s]\neq [0]$ and
$\delta([s])=1$ if $[s]=[0]$.

Since $j$ is fixed, let
$$ F =\sum_{n\in I_j} \mathbf{1}_{V_{j,n}-\sigma_{j,n}}
\quad\text{and} \quad
F_M =\sum_{n\in I_{j,M}} \mathbf{1}_{V_{j,n}-\sigma_{j,n}},$$
where the finite sets $I_{j,M}$ increase to $I_j$ as
$M\rightarrow\infty$.
We shall show that $F\in L^1(H^{\perp})$, noting that
$\text{supp}\, F_M \subseteq\text{supp}\, F \subseteq H^{\perp}$.
First,
\begin{equation}
\label{eq:3.BB}
\int_{H^{\perp}} F_M(\gamma) d\gamma
= \sum_{n\in I_{j,M}} \int_{H^{\perp}}
\mathbf{1}_{V_{j,n}-\sigma_{j,n}}(\gamma) d\gamma
= \sum_{n\in I_{j,M}} \nu(V_{j,n}) ;
\end{equation}
and it is clear that
\begin{equation}
\label{eq:3.CC}
0\leq F_M \leq F
\quad\text{and}\quad
\lim_{M\rightarrow\infty} F_M = F \quad \nu-\text{a.e.}
\end{equation}
The Beppo Levi theorem applies by \eqref{eq:3.CC},
and so, using \eqref{eq:3.BB},
$$
\int_{H^{\perp}} F(\gamma) d\gamma
= \lim_{M\rightarrow\infty} \int_{H^{\perp}} F_M(\gamma) d\gamma
= \lim_{M\rightarrow\infty} \sum_{n\in I_{j,M}} \nu(V_{j,n}) = 1 .
$$
Thus, $F\in L^1(H^{\perp})$ since $F\geq 0$.
Hence, we can apply the Lebesgue dominated convergence theorem
to the sequence $\{F_M(\cdot) \overline{ (s, \cdot ) } \}$,
which
converges $\nu$-a.e.~to $F(\cdot) \overline{ (s, \cdot ) }$
and which is bounded in absolute value by $F\in L^1(H^{\perp})$.
Therefore, \eqref{eq:3.AA} becomes
$$
\delta([s]) = \lim_{M\rightarrow\infty} \sum_{n\in I_{j,M}}
\int_{V_{j,n}-\sigma_{j,n}} \overline{ (s, \beta )} d\beta
= \lim_{M\rightarrow\infty}
\int_{H^{\perp}} F_M(\beta) \overline{ (s, \beta )} d\beta
= \int_{H^{\perp}} F(\beta) \overline{ (s, \beta )} d\beta,
$$
where we can once again consider $\nu_{H^{\perp}}$ as
the restriction of $\nu$ to $H^{\perp}$
since $H^{\perp}$ is open in $\widehat{G}$.

By the $L^1$-uniqueness theorem
for Fourier series, we must have $F=1$ on $H^{\perp}$.
Moreover, because
$F= \sum_{n\in I_j} \mathbf{1}_{V_{j,n}-\sigma_{j,n}}$,
it follows that
$\{V_{j,n}-\sigma_{j,n}\}$ is a partition
of $H^{\perp}$.  Hence, $\{V'_{j,n}-\sigma'_0\}$,
and therefore $\{V'_{j,n}\}$, are partitions
of $H^{\perp}$.
\end{proof}

\begin{remark}
a. The content of part {\em iv} is
the assertion that either the necessary or the sufficient
conditions in Theorem~\ref{thm:ONB} imply
that $\widehat{G}$ is $\sigma$-compact.  The proof
we give shows this implication when assuming
that $\{\Omega_1,\ldots,\Omega_N\}$ is a
wavelet collection of sets; the fact that the
sufficient conditions also imply $\sigma$-compactness
then follows from part {\em iii}.
However, it is easy to check
that condition {\em a} of Theorem~\ref{thm:ONB}
implies that $\widehat{G}$ is $\sigma$-compact directly,
without appealing to {\em iii}.

b. The proof that conditions {\em a} and {\em b} are
sufficient for a wavelet collection of sets
does not {\em a priori} require
$\widehat{G}$ to be $\sigma$-compact,
even though condition {\em a} easily
implies $\sigma$-compactness,
and condition {\em b} implies that
each $\nu(\Omega_j)=1$.
On the other hand, the $\sigma$-compactness
is required in the proof, in part {\em v}, that the conditions
are necessary.

c. Part {\em iv} shows that $\widehat{G}$
is $\sigma$-compact by showing that $\calD$,
and hence $\widehat{G}/H^{\perp}$, is countable.
In fact, these two conditions are equivalent;
the proof is elementary.

d. Recall that a LCAG $G$ is metrizable if and only if
$\widehat{G}$ is $\sigma$-compact.

e. Also recall, as noted before Lemma~\ref{lem:aut}, that
if there is an expansive automorphism $A\in\Aut G$,
then $\widehat{G}$ is $\sigma$-compact.
\end{remark}

\section{The Construction of Wavelet Sets}
\label{sec:thm2}

\subsection{The algorithm}
\label{ssec:const}

We are now prepared to construct wavelet collections of sets.
As before, $G$ is a LCAG with compact open subgroup $H$,
and $\calD$ is a choice of coset representatives in $\widehat{G}$
for the quotient $\widehat{G}/H^{\perp}$.
Let $A_1\in\Aut G$ be expansive with respect to $H$,
and let $\calA = \{A_1^n:n\in\ZZ\}$ be the subgroup of
$\Aut G$ generated by $A_1$.

Let $M\geq 0$ be any nonnegative
integer, and set $W=(A_1^{\ast})^M H^{\perp}$.
Note that $W\subsetneq (A_1^{\ast}) W$, because
$A_1$ is expansive.
We shall be particularly interested in the set
$({A}_1^{\ast} W) \setminus W$.
Informally speaking, $W$ is like a disk centered at the
origin, so that $({A}_1^{\ast} W) \setminus W$ is
like an annulus.  Our algorithm, based on similar
constructions for $\RR^d$ in \cite{BL99, BL01},
will translate subsets of $W$ into
$({A}_1^{\ast} W) \setminus W$
until what remains is a wavelet collection of sets.

Let $N\geq 1$ be a positive integer, which ultimately
will be the number of wavelet generators.
Let $\Omega_{1,0},\ldots , \Omega_{N,0} \subseteq W$
be measurable sets, each of
which is $(\tau,\calD)$-congruent to $H^{\perp}$,
and let 
$\widetilde{\Omega}_0 = \bigcup_{j=1}^N \Omega_{j,0}$.
We do not require the sets $\Omega_{j,0}$ to be
disjoint, but we assume that there is some integer
$\ell \geq 0$ such that
$(A^{\ast})^{-\ell} H^{\perp} \subseteq \widetilde{\Omega}_0$.
Thus, the union of all the sets $\Omega_{j,0}$ contains
a neighborhood of the origin.

For each $j\in\{1,\ldots, N\}$,
let $T_j: W \rightarrow ({A}_1^{\ast} W) \setminus W$
be a measurable {\em injective} function such that
\begin{equation}
\label{eq:Tdef}
\forall\, \gamma\in W, \quad
\quad T_j(\gamma) - \gamma + \theta_{\calD}(\gamma) \in
\calD\cap [({A}_1^{\ast} W) \setminus W] ,
\end{equation}
where $\theta_{\calD}(\gamma)$ is the unique element
of $\calD\cap (\gamma + H^{\perp})$, as in
Definition~\ref{def:DetaW}.
Thus, because $\calD\cap({A}_1^{\ast} W)$
is finite, $T_j$ slices $W$ into finitely
many measurable pieces $U_{j,n}$, each of which it translates by
some element of the form $\sigma_{j,n} - \sigma'_{j,n}$, so
that each piece lands in the ``annulus'' $({A}_1^{\ast} W) \setminus W$,
and in such a way that the translated pieces do not overlap.
Moreover, $\sigma'_{j,n}\in \calD\cap W$ and
$\sigma_{j,n}\in\calD\cap [({A}_1^{\ast} W) \setminus W]$,
with $U_{j,n}\subseteq \sigma'_{j,n} + H^{\perp}$, so that
$T_j(U_{j,n})\subseteq \sigma_{j,n} + H^{\perp}$.

Finally, we require the following compatibility condition
between the functions $T_j$ and the sets $\Omega_{j,0}$.
For any distinct $j,k\in\{1,\ldots, N\}$, assume that
either
\begin{equation}
\label{eq:Tcond1}
T_j W\cap T_k W = \emptyset,
\end{equation}
or
\begin{equation}
\label{eq:Tcond2}
T_j = T_k \quad\text{and}\quad
\Omega_{j,0}\cap\Omega_{k,0} = \emptyset .
\end{equation}

Given the initial data described above, the algorithm
proceeds as follows.
First, for each $j\in\{1,\ldots, N\}$, define
\begin{alignat}{3}
\label{eq:Lambda1}
     \Lambda'_{j,1} &= \Omega_{j,0} \cap \left[
     \bigcup_{n\geq 1} (A_1^{\ast})^{-n} \widetilde{\Omega}_{0} \right] ,
     & &
     & \Lambda''_{j,1} & = \Omega_{j,0} \cap \left[
     \bigcup_{k=1}^{j-1} \Omega_{k,0} \right] , \\
\label{eq:Lambda2}
     \Lambda_{j,1} &= \Lambda'_{j,1} \cup \Lambda''_{j,1} ,
     & & \quad\text{and}\quad
     & \Omega_{j,1} &= \left(\Omega_{j,0} \setminus \Lambda_{j,1} \right)
     \cup T_j \Lambda_{j,1}.
\end{alignat}
In addition, define
$$\widetilde{\Lambda}'_1 = \bigcup_{k=1}^N \Lambda'_{k,1}.$$
Next, for every $m\geq 1$ and $1\leq j\leq N$, define
\begin{equation}
\label{eq:alg1}
     \Lambda_{j,m+1} = \Omega_{j,m} \cap \left[
     \bigcup_{n\geq 1} (A_1^{\ast})^{-n} \widetilde{\Omega}_{m} \right]
     \quad\mbox{and}\quad
     \Omega_{j,m+1} = \left(\Omega_{j,m} \setminus \Lambda_{j,m+1} \right)
     \cup T_j \Lambda_{j,m+1},
\end{equation}
where
$$\widetilde{\Omega}_m = \bigcup_{k=1}^N \Omega_{k,m} .$$
Finally, for every $j\in\{1,\ldots, N\}$, define
\begin{equation}
\label{eq:alg4}
     \Lambda_j = \bigcup_{m\geq 1} \Lambda_{j,m} , \quad
     \Omega_j = (\Omega_{j,0} \setminus \Lambda_j) \cup T_j\Lambda_j ,
     \quad\mbox{and}\quad
     \widetilde{\Omega} = \bigcup_{j=1}^N \Omega_j .
\end{equation}
Theorem~\ref{thm:wvltset} below will show
that the $\{\Omega_1,\ldots, \Omega_N\}$ so constructed
is a wavelet collection of sets.

The inductive definitions above make sense provided that
each $\Lambda_{j,m}$ lies in $W$, the domain of $T_j$.
To see that this is indeed the case,
note that, by a simple induction on $m$, we have
\begin{equation}
\label{eq:alg2}
     \Lambda_{j,m}\subseteq \Omega_{j,0}
     \quad\mbox{and}\quad
     \Omega_{j,m}\subseteq \left[ ({A}_1^{\ast} W)
     \setminus W \right] \cup \Omega_{j,0}.
\end{equation}
Thus, $\Lambda_{j,m}\subseteq \Omega_{j,0}\subseteq W$,
as desired.
By another simple induction, we also observe that
\begin{equation}
\label{eq:Omega}
     \Omega_{j,m} =
     \left(\Omega_{j,0} \setminus \bigcup_{i=1}^m \Lambda_{j,i} \right)
     \cup T_j \left(\bigcup_{i=1}^m \Lambda_{j,i} \right).
\end{equation}

The algorithm is never vacuous, in the sense that
given $G$, $H$, $\calD$, and $A_1$,
a choice of $M$, $N$, $\{T_j\}$, and $\{\Omega_{j,0}\}$
satisfying our criteria
always exists.  For example, note that
$H^{\perp} \subsetneq A_1^{\ast} H^{\perp}$ implies
that $(A_1^{\ast} H^{\perp} ) \setminus H^{\perp}$
must contain an element $\sigma_1$ of $\calD$.
Thus, letting $\sigma'_0$ denote the unique element
of $\calD\cap H^{\perp}$, we may choose $M=0$ (so $W=H^{\perp}$),
$N=1$,
$T_1(\gamma)=\gamma-\sigma'_0 + \sigma_1$ for all $\gamma\in W$,
and $\Omega_{1,0}=H^{\perp}$.
Of course, many other choices are possible as well,
as we shall see in Section~\ref{sec3}.

\subsection{Validity of the algorithm}
\label{ssec:algpf}

Intuitively, the algorithm just described should generate a
wavelet collection of sets for the following reasons.  The
structure of $T_j$ will force each
$\{\Omega_{j,m+1}\}$ to be $(\tau,\calD)$-congruent
to $\Omega_{j,m}$ and hence, inductively, to $H^{\perp}$.
At the same time, each $\widetilde{\Omega}_m$ shares with
$\widetilde{\Omega}_0$ the property that 
$\bigcup_{n\in\ZZ} (A_1^{\ast})^n\widetilde{\Omega}_m$
covers $\widehat{G}$ up to sets of measure zero.
Meanwhile, the sets $\Lambda''_{j,1}$ guarantee that
the sets $\{\Omega_{j,1}\}_{j=1}^N$ are pairwise disjoint; and
conditions~\eqref{eq:Tcond1} and~\eqref{eq:Tcond2}
extend that property to $\{\Omega_{j,m}\}_{j=1}^N$,
for every $m\geq 1$.
More generally, the sets $\Lambda_{j,m}$ are regions
corresponding to overlaps between
different dilates $(A_1^{\ast})^n \Omega_{j,m}$.
By translating them out to $(A_1^{\ast} W)\setminus W$,
future overlaps should be smaller, and should vanish
in the limit.

Obviously, the rigorous proof requires more detail,
including the following lemma.

\begin{lemma}
\label{lem:lamblem}
Given $\{\Lambda_{j,1}\}$, $\{\Lambda_j\}$,
and $\widetilde{\Omega}$ as defined
in \eqref{eq:Lambda1}--\eqref{eq:alg4} above.
Then for every $j=1,\ldots, N$,
$$\Lambda_j \subseteq \Lambda_{j,1} \cup \left[
\bigcup_{n\geq 1}(A_1^{\ast})^{-n} \widetilde{\Omega} \right].$$
\end{lemma}

\begin{proof}
By \eqref{eq:alg4} and~\eqref{eq:Omega},
we have $\Omega_{j,m} \subseteq \Omega_{j,0} \cup T_j\Lambda_j
\subseteq \Omega_{j,0} \cup \Omega_j$,
for any $m\geq 0$ and $1\leq j\leq N$.
It follows that
$\widetilde{\Omega}_{m}\subseteq \widetilde{\Omega}_0 \cup \widetilde{\Omega}$
for any $m\geq 1$.
By \eqref{eq:alg2} and the
fact that $W \subseteq A_1^{\ast} W$, we have
$(A_1^{\ast})^{-n} \Omega_{j,m}\subseteq W$, and hence
$(T_j \Lambda_j) \cap (A_1^{\ast})^{-n} \widetilde{\Omega}_{m}=\emptyset$
for all $m\geq 0$, $n\geq 1$, and $1\leq j \leq N$.
Thus,
\begin{align*}
     \Lambda_j &= \bigcup_{m\geq 1} \Lambda_{j,m}
     = \Lambda_{j,1} \cup \bigcup_{m\geq 1} \left[
     \Omega_{j,m} \cap \left(
     \bigcup_{n\geq 1} (A_1^{\ast})^{-n} \widetilde{\Omega}_{m} \right) \right]
\\
     & \subseteq \Lambda_{j,1} \cup
     \bigcup_{m\geq 1} \left[
     (\Omega_{j,0} \cup T_j \Lambda_j) \cap \left(
     \bigcup_{n\geq 1} (A_1^{\ast})^{-n} \widetilde{\Omega}_{m} \right) \right]
\\
     & = \Lambda_{j,1} \cup
     \bigcup_{m\geq 1} \left[
     \Omega_{j,0} \cap \left(
     \bigcup_{n\geq 1} (A_1^{\ast})^{-n} \widetilde{\Omega}_{m} \right) \right]
\\
     & \subseteq \Lambda_{j,1} \cup \left[
     \Omega_{j,0} \cap \left( \bigcup_{n \geq 1}
     (A_1^{\ast})^{-n} ( \widetilde{\Omega}_0 \cup \widetilde{\Omega} )
     \right) \right]
     \subseteq \Lambda_{j,1} \cup
     \left[ \bigcup_{n\geq 1} (A_1^{\ast})^{-n} \widetilde{\Omega} \right],
\end{align*}
where the final inclusion follows because
$\Omega_{j,0} \cap \left( \bigcup_{n \geq 1}
(A_1^{\ast})^{-n} \widetilde{\Omega}_0 \right)
=\Lambda'_{j,1}\subseteq\Lambda_{j,1}$.
\end{proof}

\begin{theorem}
\label{thm:wvltset}
Let $G$ be a LCAG with
compact open subgroup $H \subseteq G$,
let $\calD\subseteq\widehat{G}$
be a choice of coset representatives in $\widehat{G}$
for $\widehat{H}=\widehat{G}/H^{\perp}$,
let $A_1\in\Aut G$ be an expansive automorphism,
and set $\calA= \{A_1^n : n\in\ZZ \}$.
For a fixed nonnegative integer $M\geq 0$
define the set $W=(A_1^*)^M H^{\perp}$;
and for a fixed positive integer $N\geq 1$
let $\{T_j: j=1,\ldots,N \}$
be a set of measurable injective functions
$$T_j : W \rightarrow (A_1^{\ast} W) \setminus W ,
\quad j=1,\ldots, N$$
satisfying \eqref{eq:Tdef}.
Consider a finite sequence
$\Omega_{1,0}, \ldots , \Omega_{N,0} \subseteq W$
of measurable subsets of $\widehat{G}$,
each of which is $(\tau,\calD)$-congruent
to $H^{\perp} \subseteq \widehat{G}$, and for which
there is a nonnegative integer $\ell\geq 0$
such that
$$(A_1^{\ast})^{-\ell} H^{\perp} \subseteq
\widetilde{\Omega}_0=\bigcup_{j=1}^N \Omega_{j,0} . $$
Assume that for every distinct pair $j,k\in\{1,\ldots, N\}$,
either condition~\eqref{eq:Tcond1} or~\eqref{eq:Tcond2}
holds.
Then the sequence $\{\Omega_1,\ldots, \Omega_N\}$ produced by the 
algorithm~\eqref{eq:Lambda1}--\eqref{eq:alg4} forms a wavelet
collection of sets.
\end{theorem}

\begin{proof}
We shall use the criterion from Theorem~\ref{thm:ONB}
to prove this theorem.

{\em i.}  It follows easily from the definition
of the maps $T_j$ that if $X$ is a measurable subset
of $W$, then $T_j(X)$ is a measurable subset of $\widehat{G}$.
The sets $\{\Omega_j\}$
are formed from the measurable sets $\{\Omega_{j,0}\}$
using only complements, countable unions, countable intersections,
and images under the maps $T_j$.
Hence, each $\Omega_j$ is measurable.

{\em ii.}
We shall now show that for any $j=1,\ldots, n$,
$\Omega_j$ is $(\tau,\calD)$-congruent to $H^{\perp}\subseteq \widehat{G}.$

Fix $j\in\{1,\ldots,N\}$.
Let $I=\calD\cap[({A}_1^{\ast} W) \setminus W]$
and $I'=\calD\cap W$, both of which are finite sets.
For each pair $(\sigma,\sigma')\in I\times I'$, define
$$
    U_{j,\sigma,\sigma'} =
    \{\gamma \in \sigma' + H^{\perp} :
    T_j(\gamma) = \gamma - \sigma' + \sigma \} .
$$
Then by definition of $T_j$, each $\gamma\in W$ lies in
exactly one $U_{j,\sigma,\sigma'}$.  Moreover,
the translated sets $T_j(U_{j,\sigma,\sigma'})$ are
pairwise disjoint, by the hypothesis that $T_j$ is injective.

For each $(\sigma,\sigma')\in I\times I'$, define
$V'_{j,\sigma,\sigma'} = (U_{j,\sigma,\sigma'}\cap\Lambda_j)$
and $V_{j,\sigma,\sigma'}=T_j(V'_{j,\sigma,\sigma'})$.
In addition, for every $\sigma'\in I'$,
define $V'_{j,\sigma'} = V_{j,\sigma'}
= (\Omega_{j,0}\setminus\Lambda_j) \cap (\sigma' + H^{\perp})$.
Using the finite set $I' \cup (I\times I')$ as an indexing
set, we are now prepared to exhibit the 
$(\tau,\calD)$-congruence.

Because $\Lambda_j\subseteq\Omega_{j,0}$,
$\{V'_{j,\sigma,\sigma'} : (\sigma,\sigma')\in I\times I'\} \cup
\{V'_{j,\sigma'} : \sigma'\in I'\}$ 
tiles $\Omega_{j,0}$.  By definition of $\Omega_j$,
$\{V_{j,\sigma,\sigma'} : (\sigma,\sigma')\in I\times I'\} \cup
\{V_{j,\sigma'} : \sigma' \in I'\}$
tiles $\Omega_j$.
Moreover,
$V_{j,\sigma,\sigma'} -\sigma = V'_{j,\sigma,\sigma'} -\sigma'$
and
$V_{j,\sigma'} = V'_{j,\sigma'}$.
Therefore, by definition of $T_j$,
$\Omega_j$ is $(\tau,\calD)$-congruent to $\Omega_{j,0}$.
Since $\Omega_{j,0}$ is $(\tau,\calD)$-congruent to $H^{\perp}$,
it follows that
$\Omega_j$ is $(\tau,\calD)$-congruent to $H^{\perp}$.

{\em iii.}
Next, we shall show that for any $j,k\in\{1,\ldots, N\}$
and for any {\em distinct} $A,B\in\calA$,
the sets $A^{\ast}\Omega_j$ and $B^{\ast}\Omega_k$ are disjoint.
Because $\calA=\{A_1^n : n\in\ZZ\}$, it suffices to show that
$$ 
   \Omega_k \cap (A_1^{\ast})^{n} \Omega_j =\emptyset
   \quad\text{for all } j,k\in\{1,\ldots, N\}
   \text{ and } n\geq 1.
$$

Fix $j,k\in\{1,\ldots, N\}$ and $n\geq 1$, and suppose that there is some
$\gamma\in \Omega_k \cap (A_1^{\ast})^{n}\Omega_j .$
Thus, $\gamma\in \Omega_k$ and
$(A_1^{\ast})^{-n} \gamma \in \Omega_j$.

If $\gamma\in \Omega_{k,0}$, then because
$\Omega_{k,0}\subseteq W \subseteq (A_1^{\ast})^{n} W$, we have
$(A_1^{\ast})^{-n} \gamma \in W\cap\Omega_j \subseteq\Omega_{j,0}$,
by \eqref{eq:alg4}
and the fact that the image of $T_j$ is disjoint from $W$.
Thus,
$$
(A_1^{\ast})^{-n} \gamma \in
\Omega_{j,0}\cap (A_1^{\ast})^{-n} \Omega_{k,0}
\subseteq \Omega_{j,0}\cap (A_1^{\ast})^{-n} \widetilde{\Omega}_0
\subseteq \Lambda_{j,1}
\subseteq \Lambda_j.
$$
Then $(A_1^{\ast})^{-n}\gamma\not\in \Omega_j$,
which is a contradiction.

On the other hand, if $\gamma\not\in\Omega_{k,0}$, then
$\gamma\in T_k\Lambda_k$ by \eqref{eq:alg4},
so that $\gamma\in T_k\Lambda_{k,m}$
for some $m\geq 1$.
Note that $\gamma\in A_1^{\ast} W$, so
that $(A_1^{\ast})^{-n} \gamma\in W$, and therefore
$$(A_1^{\ast})^{-n}\gamma
\in W \cap \Omega_j
\subseteq \Omega_{j,0} \setminus \Lambda_j
\subseteq \Omega_{j,0} \setminus \left( \bigcup_{i=1}^m \Lambda_{j,i} \right)
\subseteq  \Omega_{j,m},$$
by \eqref{eq:alg4} and \eqref{eq:Omega}.
Also note that
$\gamma\in T_k \Lambda_{k,m} \subseteq\Omega_{k,m}
\subseteq\widetilde{\Omega}_m$.
Thus, we have
$\gamma\in\widetilde{\Omega}_m$ and
$(A_1^{\ast})^{-n}\gamma \in \Omega_{j,m}$; hence,
$$
     (A_1^{\ast})^{-n}\gamma \in \Omega_{j,m}
     \cap (A_1^{\ast})^{-n} \widetilde{\Omega}_m
     \subseteq\Lambda_{j,m+1}\subseteq\Lambda_j,
$$
contradicting the assumption that $(A_1^{\ast})^{-n}\gamma\in\Omega_j$.

{\em iv.}
Now we shall show that for any distinct pairs
$(A,j),(B,k)\in \calA \times \{1,\ldots, N\}$,
the sets $A^{\ast}\Omega_j$ and $B^{\ast}\Omega_k$ are disjoint.
In light of {\em iii}, we may assume that $A=B$, and therefore
it suffices to show that for any distinct $j,k\in\{1,\ldots, N\}$,
the sets $\Omega_j$ and $\Omega_k$ are disjoint.

Without loss, assume that $k<j$.  If $T_j=T_k$,
then by hypothesis~\eqref{eq:Tcond2},
$\Omega_{j,0}\cap\Omega_{k,0}=\emptyset$,
and therefore $\Lambda_j\cap\Lambda_k=\emptyset$,
by \eqref{eq:alg2}.  Because $T_j=T_k$ is injective, it
follows that $T_j\Lambda_j \cap T_k\Lambda_k =\emptyset$.
Moreover, the image of $T_j$ does not intersect
$\Omega_{j,0},\Omega_{k,0}\subseteq W$.
It follows that $\Omega_j \cap \Omega_k = \emptyset$.

On the other hand, if $T_j\neq T_k$,
then by hypothesis~\eqref{eq:Tcond1},
$T_j W \cap T_k W =\emptyset$,
so that
$T_j \Lambda_j \cap T_k \Lambda_k =\emptyset$.
Moreover, 
$\Omega_{j,0}\cap\Omega_{k,0}
\subseteq \Lambda''_{j,1}\subseteq \Lambda_j$
by~\eqref{eq:Lambda1} and~\eqref{eq:Lambda2},
so that
$$(\Omega_{j,0}\setminus \Lambda_j)
\cap(\Omega_{k,0}\setminus \Lambda_k) =\emptyset .$$
Again by the fact that 
the images of $T_j$ and $T_k$ do not intersect
$\Omega_{j,0},\Omega_{k,0}\subseteq W$,
it follows that $\Omega_j \cap \Omega_k = \emptyset$.

{\em v.}
Finally, we need to show that
$\{ A^{\ast}\Omega_j : A\in\calA,\, 1\leq j \leq N \}$
covers $\widehat{G}$ up
to sets of measure zero.  Let
$Y=\bigcap_{n\in\ZZ} (A_1^{\ast})^{n} H^{\perp}$.
Because $\nu((A_1^{\ast})^{n} H^{\perp} )=|A_1|^n$, we have
$$\nu(Y)\leq\inf\{|A_1|^n : n\in\ZZ\} = 0 .$$
Therefore, it will suffice to show that
$$ 
    \widehat{G}\setminus Y \subseteq
    \bigcup_{A\in\calA} A^{\ast}\widetilde{\Omega} .
$$

Given $\gamma\in \widehat{G}\setminus Y$, let
$J_\gamma = \{n\in\ZZ : (A_1^{\ast})^{-n}\gamma\in \widetilde{\Omega}_0\}.$
We shall need a minimal element
of $J_{\gamma}$, and so we must show that $J_{\gamma}$ is
nonempty and bounded below.

By the hypothesis
that $(A_1^{\ast})^{-\ell} H^{\perp} \subseteq \widetilde{\Omega}_0$
and because $A_1$ is expansive,
$$
\widehat{G} =
\bigcup_{n\geq 0}  (A_1^{\ast})^{n-\ell} H^{\perp} \subseteq
\bigcup_{n\geq 0}  (A_1^{\ast})^{n} \widetilde{\Omega}_0,
$$
and therefore $J_{\gamma}$ is nonempty.
On the other hand, suppose that for every integer $n\in\ZZ$,
there is some integer $i_n\leq n$ such that $i_n\in J_{\gamma}$, that
is, $(A_1^{\ast})^{-i_n} \gamma \in \widetilde{\Omega}_0$.
Then because $\widetilde{\Omega}_0\subseteq (A_1^{\ast})^M H^{\perp}$
and because $A_1$ is expansive, we have
$$\gamma \in (A_1^{\ast})^{i_n} \widetilde{\Omega}_0
\subseteq (A_1^{\ast})^{i_n + M} H^{\perp}
\subseteq (A_1^{\ast})^{n + M} H^{\perp}.$$
Thus, $\gamma\in \bigcap_{n\in\ZZ} (A_1^{\ast})^{n + M} H^{\perp} = Y$,
which is a contradiction.  It follows that $J_{\gamma}$
is bounded below.
Hence, it makes sense to define $i=\min J_\gamma$.

We claim
that $(A_1^{\ast})^{-i}\gamma \in \widetilde{\Omega}_0
\setminus \widetilde{\Lambda}'_1$.
Indeed, we have $(A_1^{\ast})^{-i}\gamma \in \widetilde{\Omega}_0$ because
$i\in J_\gamma$.  On the other hand, if
$(A_1^{\ast})^{-i}\gamma \in \widetilde{\Lambda}'_1,$ then by definition
of $\widetilde{\Lambda}'_1$, we have
$(A_1^{\ast})^{-i} \gamma \in (A_1^{\ast})^{-n}\widetilde{\Omega}_0$ for
some $n\geq 1$.  Therefore,
$(A_1^{\ast})^{-(i-n)}\gamma \in \widetilde{\Omega}_0$,
so that $i-n\in J_{\gamma}$,
which contradicts the minimality of $i$ and proves our claim.

Let $j\in\{1,\ldots, N\}$ be the smallest index such
that $(A_1^{\ast})^{-i} \gamma\in\Omega_{j,0}$; note
that $j$ exists because $(A_1^{\ast})^{-i} \gamma\in\widetilde{\Omega}_0$.

If $(A_1^{\ast})^{-i} \gamma\in \Omega_j$, then
$\gamma\in (A_1^{\ast})^{i}\Omega_j$, and we are done.

Otherwise, $(A_1^{\ast})^{-i} \gamma\in \Lambda_j$.
However, $(A_1^{\ast})^{-i} \gamma\not\in \Lambda'_{j,1}$,
by the claim above; and
$(A_1^{\ast})^{-i} \gamma\not\in \Lambda''_{j,1}$ by
our choice of $j$.  Thus,
$(A_1^{\ast})^{-i} \gamma \in
\Lambda_j \setminus \Lambda_{j,1}$.
Therefore, by Lemma~\ref{lem:lamblem}, there is some $n\geq 1$
for which
$(A_1^{\ast})^{-i} \gamma \in (A_1^{\ast})^{-n} \widetilde{\Omega}$.
Hence, $\gamma\in (A_1^{\ast})^{i-n} \Omega$,
and the theorem is proved.
\end{proof}

\section{Examples of Sets of Wavelet Generators}\label{sec3}
\renewcommand{\theequation}{\thesection.\arabic{equation}}
\setcounter{equation}{0}

Theorems~\ref{thm:ONB} and~\ref{thm:wvltset} together provide
an algorithm for generating many wavelet orthonormal bases for
$L^2(G)$ for a group $G$ with the properties specified in
Theorem~\ref{thm:wvltset}.
In this section, we shall examine several examples
of such wavelets.
Other examples can be found in \cite{BenR03}.

\subsection{Haar/Shannon wavelets on $G$}
\label{ssec:Haar}
Let $G$ be a LCAG with compact open subgroup $H$,
let $\calD$ be a choice set of coset representatives
in $\widehat{G}$ for $\widehat{G} / H^{\perp}$,
and let $A_1$ be an expansive automorphism of $G$.

\begin{figure}
\scalebox{.8}{\includegraphics{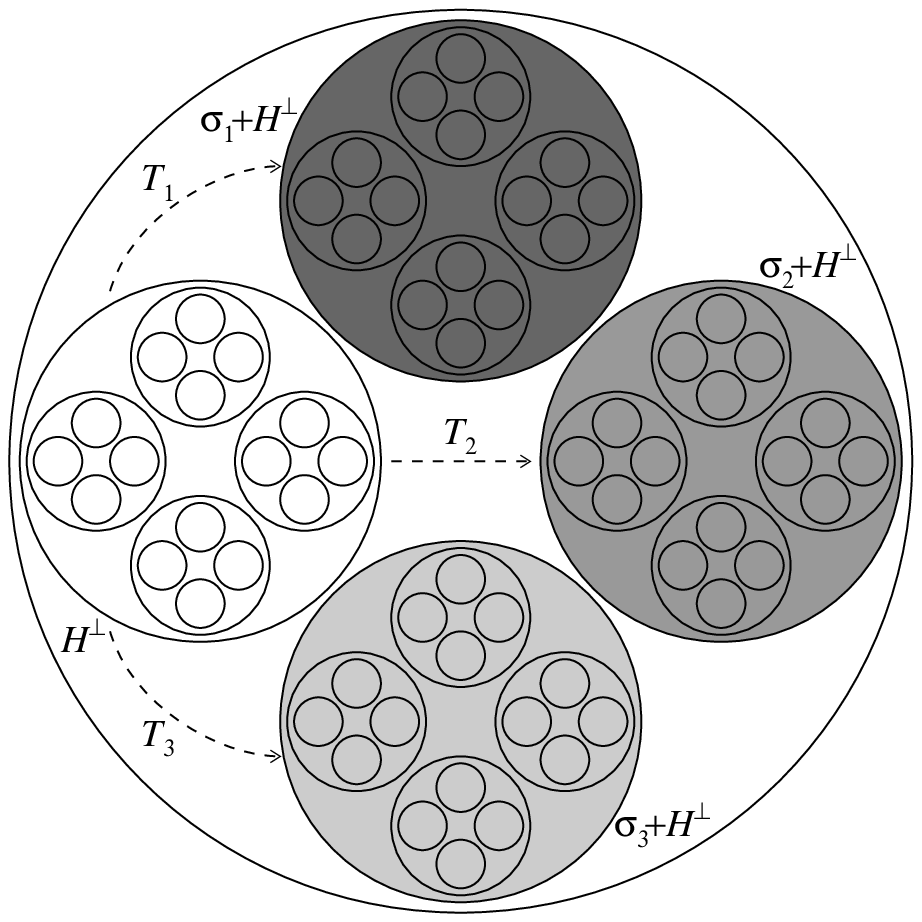}
\includegraphics{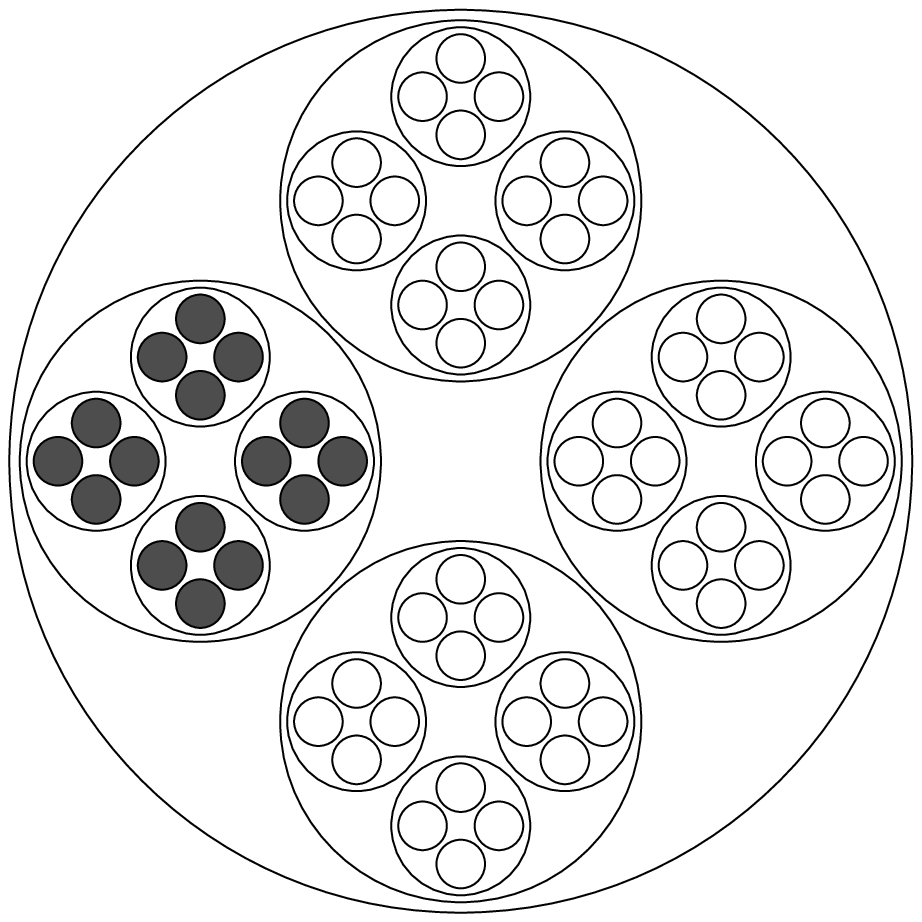}}
\caption{The maps $T_j$ and the sets $\Omega_{1,0}=\Omega_{2,0}=\Omega_{3,0}$
of Example~\ref{ssec:Haar}, for $|A_1|=4$.}
\label{fig:Tmap5.1}
\end{figure}

Take $M=0$, so that $W=H^{\perp}$,
set $N=|A_1| - 1\geq 1$, and let $\sigma_1,\ldots,\sigma_N$
be the $N$ elements of
$\calD \cap [(A_1^{\ast}H^{\perp})\setminus H^{\perp}]$.
For each $j=1,\ldots, N$,
define $T_j(\gamma)=\gamma-\sigma'_0 + \sigma_j$,
where $\sigma'_0$ denotes the unique element
of $\calD\cap H^{\perp}$,
and define $\Omega_{j,0}= H^{\perp}$.
Note that $\{H^{\perp}, T_1 H^{\perp},\ldots, T_N H^{\perp} \}$
is a set of $|A_1|=N+1$ compact open sets which together
tile $A_1^{\ast} H^{\perp}$.
See Figure~\ref{fig:Tmap5.1} for a diagram
of $\{T_j\}$ and $\Omega_{j,0}$ ($j=1,2,3$) in the case that
$|A_1|=4$.

\begin{figure}
\scalebox{.8}{\includegraphics{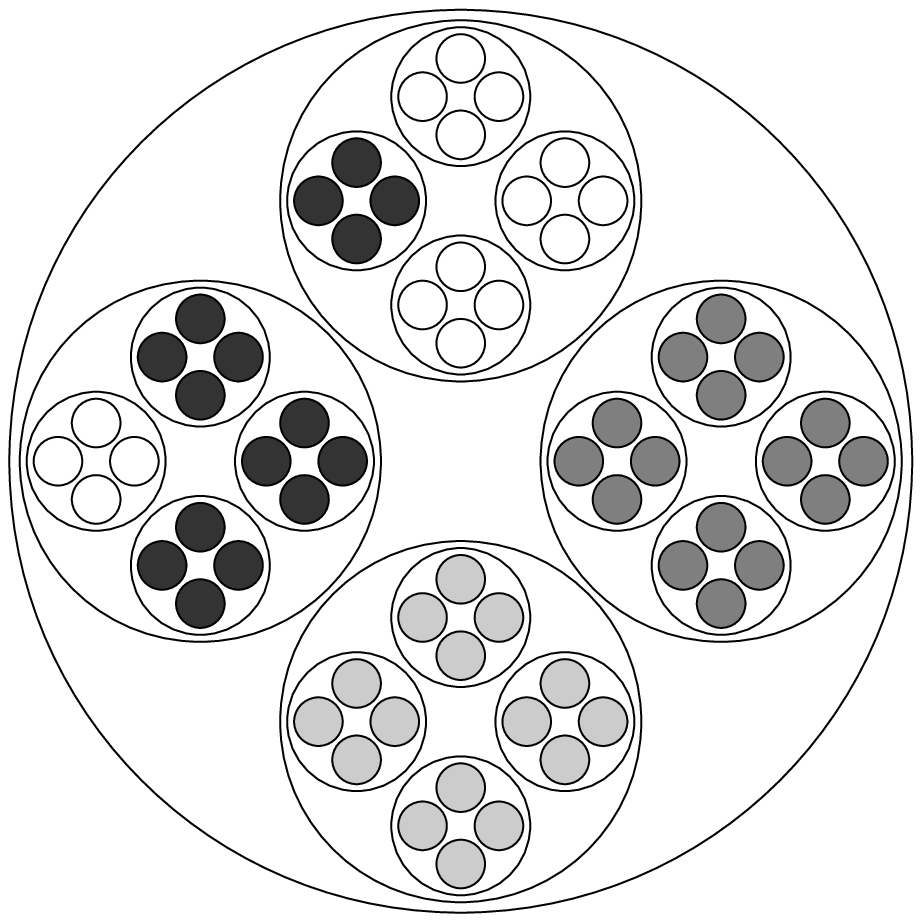}
\includegraphics{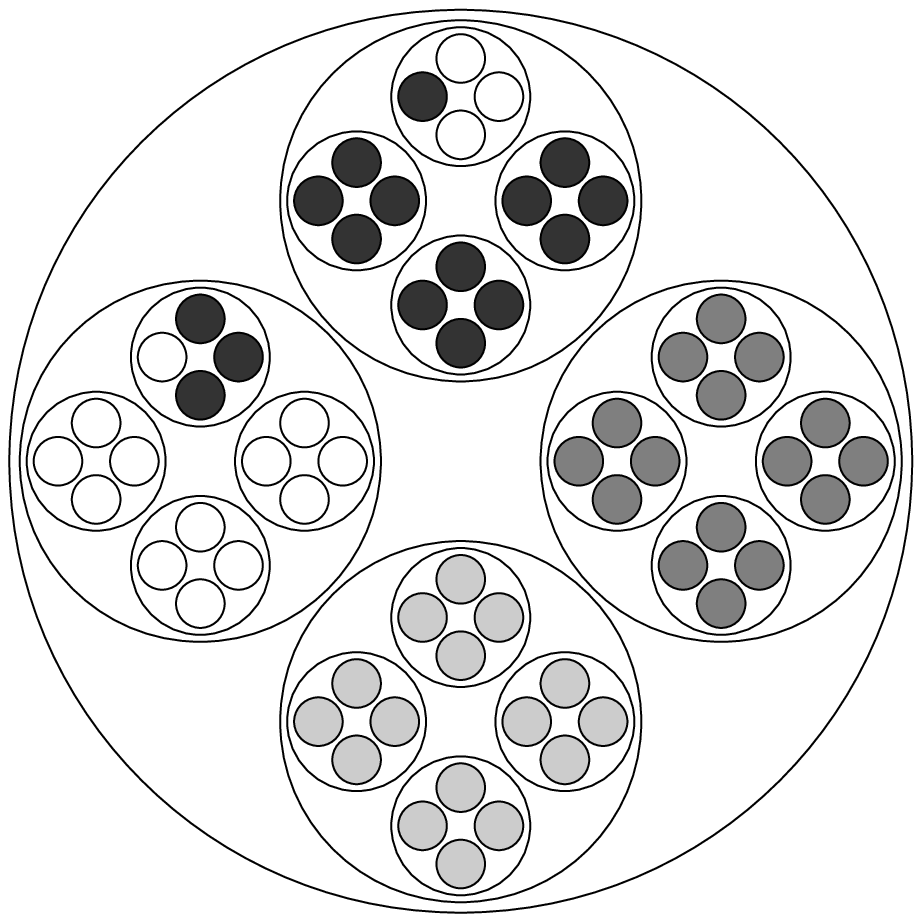}}
\caption{The sets $\Omega_{j,m}$ ($j=1,2,3$, $m=1,2$)
of Example~\ref{ssec:Haar}, for $|A_1|=4$.}
\label{fig:Omega5.1}
\end{figure}

When we apply the algorithm from Section~\ref{sec:thm2}
to this data, we obtain, up to sets of measure zero,
\begin{equation}
\label{eq:Haardef}
\Omega_j = \sigma_j - \sigma'_0 + H^{\perp}
= \sigma_j + H^{\perp} \quad \text{for all } j\in\{1,\ldots, N\}.
\end{equation}
To see this, first note that for $j\geq 2$, we have
$\Lambda''_{j,1} = H^{\perp}$, so that $\Lambda_j= H^{\perp}$,
and hence $\Omega_j=T_j H^{\perp}$, as claimed; see Figure \ref{fig:Omega5.1}.

Meanwhile, $\Lambda_{1,1}=(A_1^{\ast})^{-1}H^{\perp}$,
and the observations about $\Lambda_j$ for $j\geq 2$
imply that
$$\bigcup_{j=2}^N (A_1^{\ast})^{-1}T_j H^{\perp} \subseteq\Lambda_{1,2}.$$
Because
$\{H^{\perp}, T_1 H^{\perp},\ldots, T_N H^{\perp} \}$
tiles $A_1^{\ast} H^{\perp}$, it follows that
$$H^{\perp}\cap \Omega_{1,2} \subseteq
\gamma_2 + (A_1^{\ast})^{-1} H^{\perp} ,
\quad\text{where } \gamma_2= (A_1^{\ast})^{-1} T_1(0) .$$
By induction on $m$, we obtain
$$H^{\perp}\cap \Omega_{1,m} \subseteq
\gamma_m + (A_1^{\ast})^{-(m-1)} H^{\perp} ,
\quad\text{where } \gamma_m= (A_1^{\ast})^{-1} T_1(\gamma_{m-1}) .$$
See Figure~\ref{fig:Omega5.1} for a diagram
of $\Omega_{j,m}$ ($j=1,2,3$, $m=1,2$) in the case that
$|A_1|=4$.

Thus, $\nu(\Lambda_1) = 1-\inf\{|A_1|^{-(m-1)} : m\geq 2\} = 1$,
so that $\Lambda_1=H^{\perp}$, up to sets of measure zero.
The characterization of $\{\Omega_j : j=1,\ldots, N\}$
in \eqref{eq:Haardef} then follows.
Note that by that characterization, the set $\calD$
is in this case irrelevant to the ultimate set
of wavelet generators.
The wavelets $\{\psi_j : j=1,\ldots, N\}$ themselves are
described by the following proposition.

\begin{proposition}
\label{prop:haarshan}
Let $G$ be a LCAG with
compact open subgroup $H \subseteq G$,
let $\calD\subseteq\widehat{G}$
be a choice of coset representatives in $\widehat{G}$
for $\widehat{G}/H^{\perp}$,
let $A_1\in\Aut G$ be an expansive automorphism,
and set $N=|A_1|-1\geq 1$.
Denote by 
$\{\sigma_j : j=1,\ldots, N\}$
the $N$ elements of $\calD \cap[(A_1^* H^{\perp})\setminus H^{\perp}]$,
write $\Omega_j=\sigma_j +H^{\perp}$ for $j=1,\ldots, N$,
and define $\psi_j\in L^2(G)$ by $\widehat{\psi}_j = \mathbf{1}_{\Omega_j}$.
Then $\{\psi_1, \ldots, \psi_N\}$ is a set of wavelet generators, and,
for all $j=1,\ldots, N$ and $x\in G$,
$$\psi_j(x) = (x,\sigma_j)\mathbf{1}_H(x) .$$
Moreover, $\psi_j$ is constant on every open set of the
form $c + A_1^{-1} H$.
\end{proposition}

\begin{proof}
By the preceding argument, $\{\Omega_1, \ldots, \Omega_N\}$
is the output of algorithm \ref{eq:Lambda1} -- \ref{eq:alg4}.
By Theorem \ref{thm:wvltset}, 
$\{\psi_1, \ldots, \psi_N\}$ is a set of wavelet generators. 
We compute
\begin{equation}
\label{eq:haarcomp}
\psi_j(x)
= \int_{\widehat{G}} \mathbf{1}_{\sigma_j + H^{\perp}} (\gamma)
(x, \gamma) d\gamma
= \int_{\sigma_j + H^{\perp}} (x, \gamma) d\gamma
= (x,\sigma_j) \int_{H^{\perp}} (x, \beta) d\beta,
\end{equation}
by the change of variables $\gamma=\sigma_j + \beta$.
If $x\in H$, then $(x,\beta)=1$ for all $\beta\in H^{\perp}$,
so that the right side of \eqref{eq:haarcomp}
is $(x,\sigma_j)$, as desired.
On the other hand, if $x\not\in H$, then the integral
is zero; see the discussion following \eqref{eq:ut}.

Finally, let $c\in G$, and pick any $x,y\in c+A_1^{-1} H$.
Then $y-x\in A_1^{-1} H$; moreover, $x$ lies in $H$
if and only if $y$ does also.  Thus,
$$\psi_j(y) = (y,\sigma_j)\mathbf{1}_H(y)
= (y-x,\sigma_j) (x,\sigma_j)\mathbf{1}_H(x)
= (y-x,\sigma_j) \psi_j(x).$$
However, $(y-x,\sigma_j) = (A_1(y-x), (A_1^{\ast})^{-1}\sigma_j)$,
and by hypothesis, $A_1(y-x)\in H$ and
$(A_1^{\ast})^{-1}\sigma_j\in H^{\perp}$.
Thus, $\psi_j(y)=\psi_j(x)$, and
$\psi_j$ is constant on $c + A_1^{-1} H$, as claimed.
\end{proof}

In \cite{BenR03}, it was observed that the
wavelets $\{\psi_1,\ldots,\psi_N\}$ defined
by $\psi_j=\check{\mathbf{1}}_{\Omega_j}$
are Shannon wavelets in the sense that each
$\Omega_j$ is a simple translation of
the original ``fundamental domain'' $H^{\perp}$.
On the other hand, by Proposition~\ref{prop:haarshan},
each $\psi_j$ is a step function
with compact support, analogous to the
usual Haar wavelet.  (The functions $\psi_j$ are
also continuous, which is possible because $G$
is a totally disconnected topological space.)
Thus, although the Haar and Shannon
wavelets lie at opposite extremes in the family
of wavelets over $\RR^d$, their analogues
over $G$ are one and the same.

As a special case, let $p$ be a prime number, and
consider $G=\Qp$, as
in Example~\ref{ex:Qp}.  Let $A_1$ be
the multiplication-by-$1/p$ automorphism,
which is expansive.
Let $\calD$ be any set of coset representatives in $\widehat{G}$
for $\widehat{G} / H^{\perp} = \Qp/\Zp$.
For example, we could
use Kozyrev's set
$\calD_{\text{Koz}}=\{\frac{m}{p^e} : e\geq 1, \, 0\leq m \leq p^e - 1 \}$,
which is the set generated by Lemma~\ref{lem:Ddef}
using $\rho_0=0, \rho_1=1,\ldots,\rho_{p-1}=p-1$;
but, as observed above, the Shannon wavelets produced in the
end will be the same regardless of $\calD$.

By \eqref{eq:Haardef}, the wavelet sets generated
by the algorithm are
$$\Omega_j
= j/p + \Zp \quad \text{for } j\in\{1,\ldots, p-1\}.$$
As observed in \cite{BenR03}, one can compute that the
wavelets $\psi_j=\check{\mathbf{1}}_{\Omega_j}$
are precisely the Haar wavelets
constructed by Kozyrev \cite[Theorem 2]{Koz02}.  (On the other hand,
our translation operators are different from his,
so that the orthonormal basis is different, if only
in that some of the elements of ours are scalar
multiples of his.)

Similarly, if we let $G=\FF_2((t))$, as in Example~\ref{ex:Fpt},
then the data above gives
the (single) wavelet set $\Omega_1= 1/t + \FF_2[[t]]$,
and the corresponding wavelet is precisely the Haar wavelet
computed by Lang \cite{Lan96}.

However, even for the Haar/Shannon wavelets, our construction
is broader than those of Kozyrev and Lang, because
it applies to any group $G$ satisfying our hypotheses
(such as finite products or
finite extensions of $\Qp$ or of $\Fp((t))$),
and because it works for any given expansive automorphism.

For example, let $G=\QQ_2 \times \QQ_3$, and let $A_1$
be multiplication-by-$(1/4, 1/3)$.  Then $A_1$ is
an expansive automorphism with respect to the compact
open subgroup $H=\ZZ_2\times\ZZ_3$, and $|A_1|=12$.
The Haar/Shannon
wavelet sets in $\widehat{G}= \QQ_2 \times \QQ_3$ are
the eleven sets of the form
$$(j/4 + \ZZ_2) \times (k/3 + \ZZ_3),$$
for $0\leq j\leq 3$ and $0\leq k\leq 2$, with $j$ and $k$
not both zero.

Alternatively, let $G=\QQ_2(\sqrt{2})$, the field extension
of $\QQ_2$ obtained by adjoining $\sqrt{2}$, and let
$H=\ZZ_2[\sqrt{2}]$ be the ring of integers in this field.
Let $A_1$
be multiplication-by-$2$, which is an expansive automorphism,
this time with $|A_1|=4$.  We have $\widehat{G} = G$,
with duality given by $(x,\gamma)=\exp(2\pi i\text{Tr}(x \gamma))$,
where $\text{Tr}:\QQ_2(\sqrt{2})\rightarrow \QQ_2$ is the trace map
$\text{Tr}(y+z\sqrt{2}) = 2y$ for $y,z\in\QQ_2$.  Note that
$H^{\perp}=(2^{-3/2})\ZZ_2[\sqrt{2}]$ under this duality.
The Haar/Shannon
wavelet sets in $\widehat{G}$ are the three sets
$$ 2^{-5/2} + (2^{-3/2})\ZZ_2[\sqrt{2}], \quad
2^{-5/2} + 2^{-2} + (2^{-3/2})\ZZ_2[\sqrt{2}], \quad \text{and} \quad
2^{-2} + (2^{-3/2})\ZZ_2[\sqrt{2}].$$

\subsection{Single wavelets on $G$}
\label{ssec:single}

Let $G$ be a LCAG with compact open subgroup $H$,
let $\calD$ be a choice of coset representatives in $\widehat{G}$
for $\widehat{G} / H^{\perp}$, and
let $A_1$ be an expansive automorphism of $G$.

\begin{figure}
\scalebox{.8}{\includegraphics{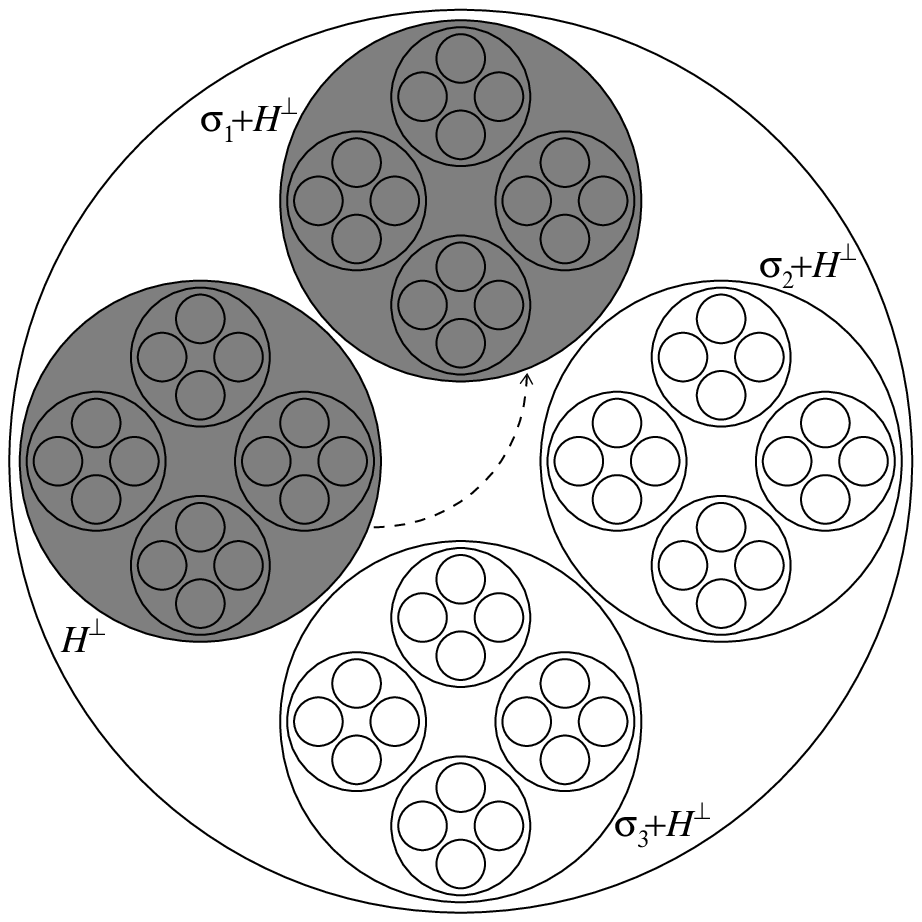}
\includegraphics{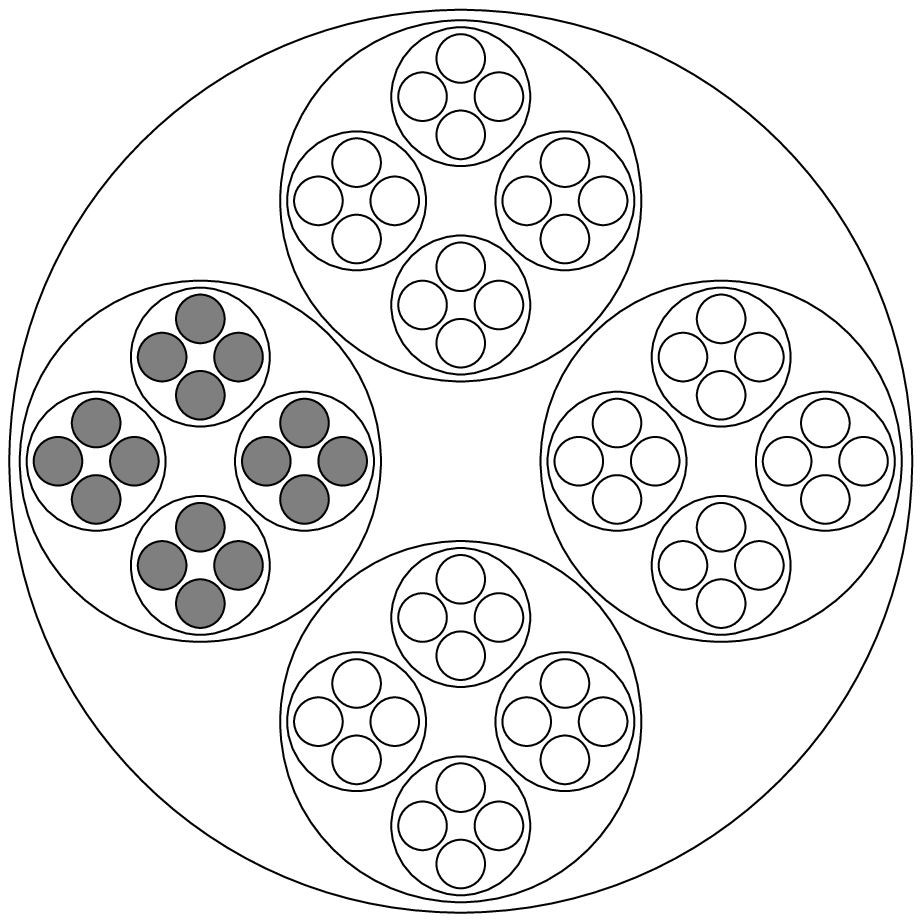}}
\caption{The map $T_1$ and set $\Omega_{1,0}$
of Example~\ref{ssec:single}, for $|A_1|=4$.}
\label{fig:Tmap5.2}
\end{figure}

Take $M=0$, so that $W=H^{\perp}$,
set $N=1$, and let $\sigma_1$
be any one of the $|A_1|-1$ elements of
$\calD\cap [(A_1^{\ast}H^{\perp})\setminus H^{\perp}]$.
Define $T_1(\gamma)=\gamma-\sigma'_0 + \sigma_1$,
where $\sigma'_0$ denotes the unique element
of $\calD\cap H^{\perp}$,
and define $\Omega_{1,0}= H^{\perp}$.
See Figure~\ref{fig:Tmap5.2} for a diagram
of $T_1$ and $\Omega_{1,0}$ in the case that
$|A_1|=4$.

\begin{figure}
\scalebox{.8}{\includegraphics{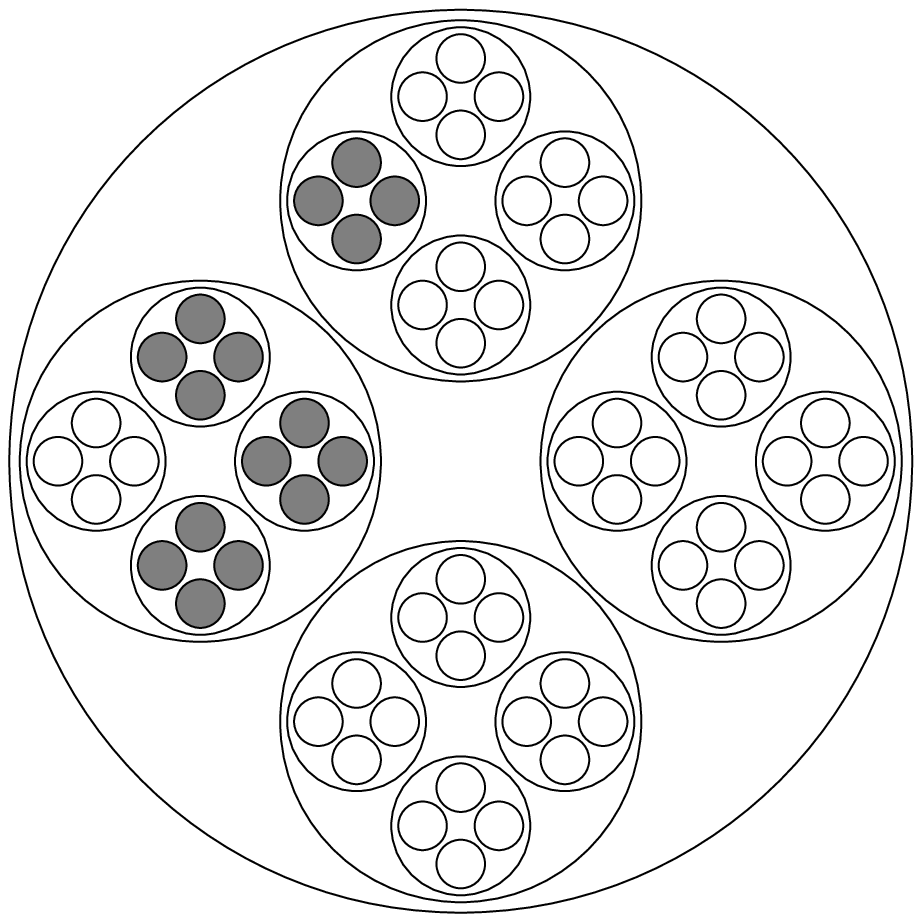}
\includegraphics{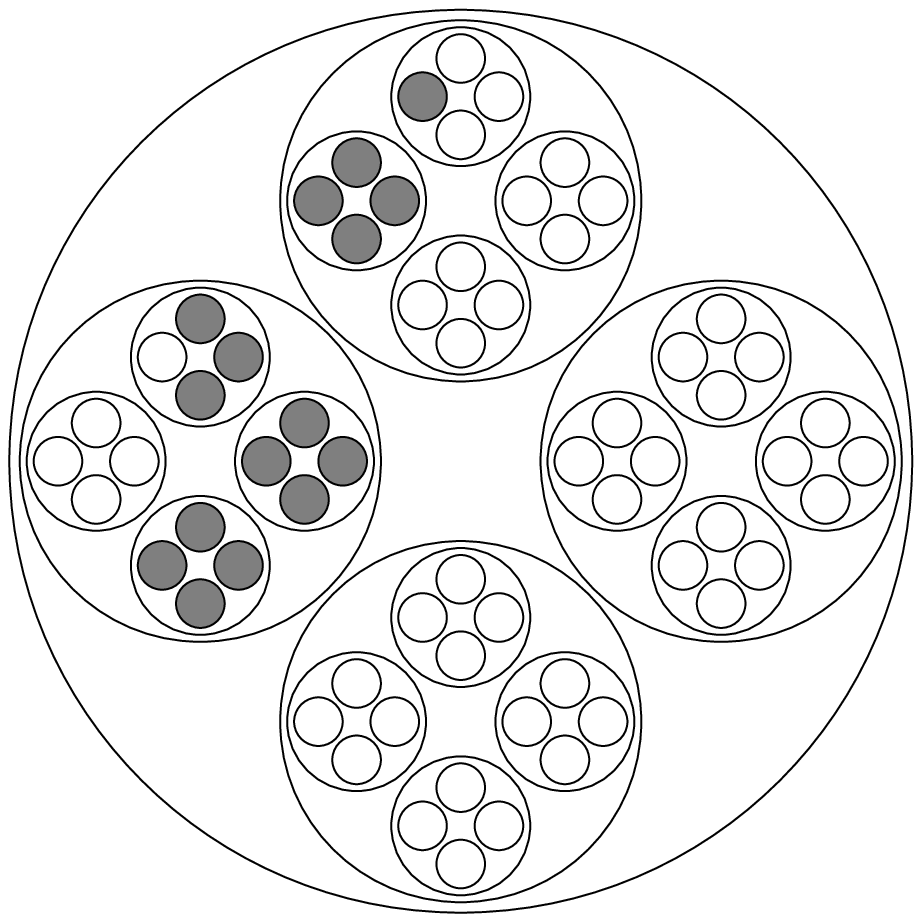}}
\caption{$\Omega_{1,1}$ and $\Omega_{1,2}$ of Example~\ref{ssec:single}.}
\label{fig:Omega5.2}
\end{figure}

A simple induction shows that for every $m\geq 1$,
$$\Lambda_{1,m} =
\left[ \sum_{i=1}^{m-1} (A_1^{\ast})^{-i} (\sigma_1 - \sigma'_0) \right]
+ (A_1^{\ast})^{-m} H^{\perp} .$$
If we think of $H^{\perp}$ as a ``disk'' of measure $1$,
then each $\Lambda_{1,m}$ is a ``disk'' of measure
$|A_1|^{-m}$.  Moreover, these disks $\Lambda_{1,m}$, $m\geq 1$,
are pairwise disjoint, and they approach the point
$$c=\sum_{i=1}^{\infty} (A_1^{\ast})^{-i} (\sigma_1 - \sigma'_0)
\in H^{\perp}.$$
See Figure~\ref{fig:Omega5.2}
for some of the resulting sets $\Omega_{1,m}$.
Thus, $\Lambda_1$ is a countable union
of such disks, and $\nu(\Lambda_1) = (|A_1| - 1)^{-1}$.
The wavelet set is then
$$\Omega_1 = ( H^{\perp}\setminus\Lambda_1 )
\cup ( \sigma_1 - \sigma'_0 + \Lambda_1 ) .$$

For example, if $G=\QQ_5$ and $A_1$ is multiplication-by-$1/5$,
and if $\calD$ contains $\sigma'_0=0$ and $\sigma_1=1/5$,
then we have
$$\Lambda_{1,m} = (1 + 5 + 5^2 + \ldots + 5^{m-2}) + 5^m \ZZ_5
= -\frac{1}{4} - 5^{m-1} + 5^m \ZZ_5 ;$$
and therefore, since $-1/4 + 1/5=1/20$, we have the wavelet set
$$\Omega_1 =
\left[ \ZZ_5 \setminus \bigcup_{m=1}^\infty 
\left( -\frac{1}{4} - 5^{m-1} + 5^m \ZZ_5 \right) \right]
\cup
\bigcup_{m=1}^\infty 
\left( -\frac{1}{20} - 5^{m-1} + 5^m \ZZ_5 \right) .$$

\subsection{Another single wavelet on $\QQ_3$}
\label{ssec:Q3wave}

Let $G=\QQ_3$, with compact open subgroup $H=\ZZ_3$,
and let $A_1$ be multiplication-by-$1/3$, so that $A_1$ is
expansive, with $|A_1|=3$.
Identify $\widehat{G}$ as $\QQ_3$ and $H^{\perp}$ as $\ZZ_3$,
as in Example~\ref{ex:Qp}.
Let $\calD$ be a set of coset representatives in $\widehat{G}$
for $\widehat{G} / H^{\perp}$ including
$\sigma'_0=0$, $\sigma_1=1/3$, and $\sigma_2 = 2/3$.

\begin{figure}
\includegraphics{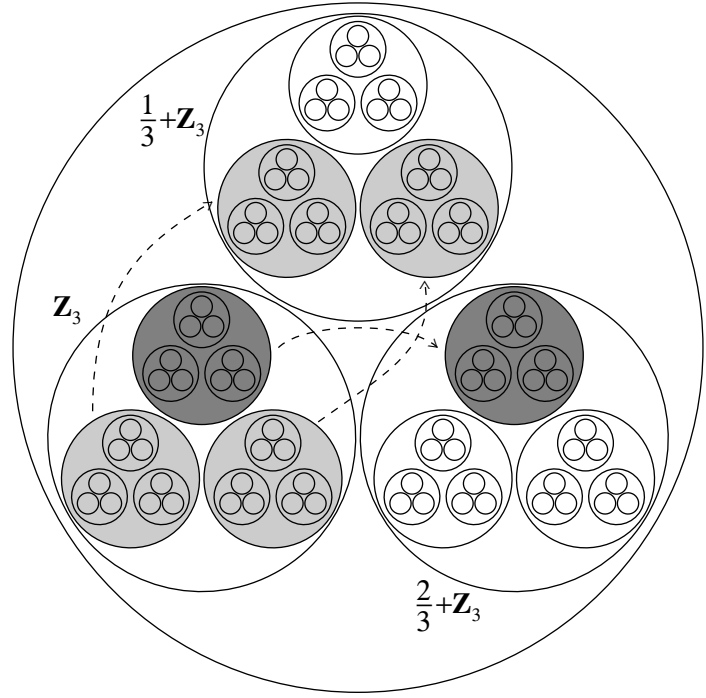}
\caption{The map $T_1$ of Example~\ref{ssec:Q3wave}.}
\label{fig:Tmap5.3}
\end{figure}

\begin{figure}
\scalebox{.8}{\includegraphics{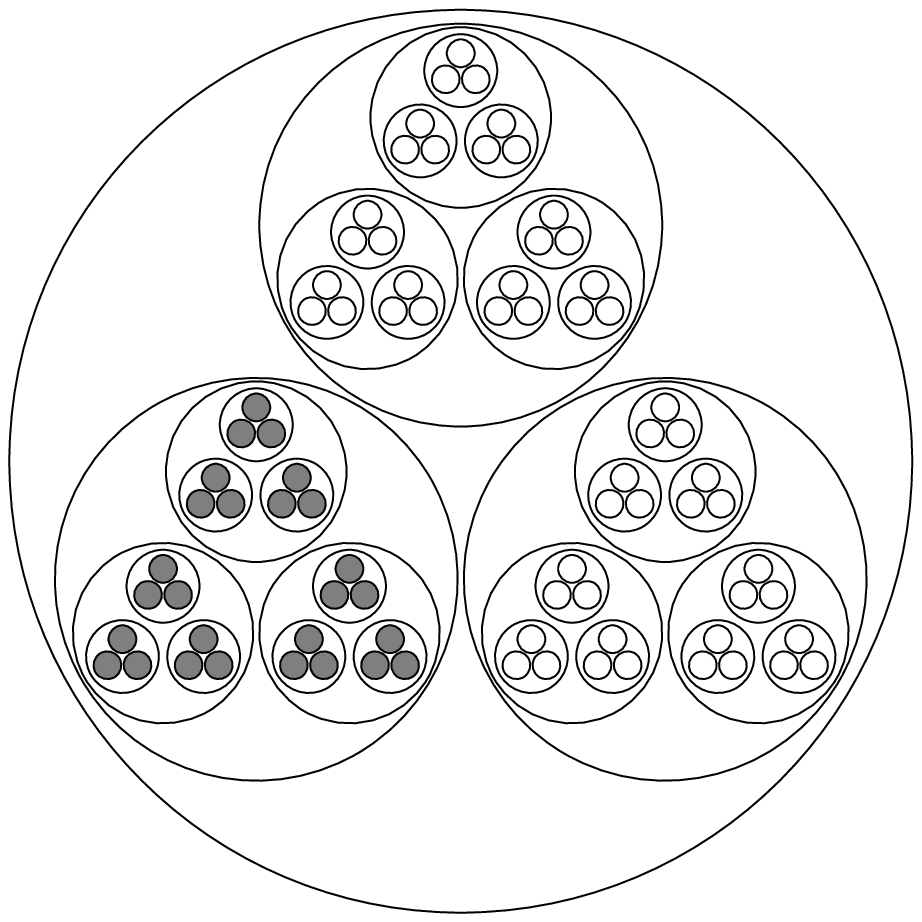}
\includegraphics{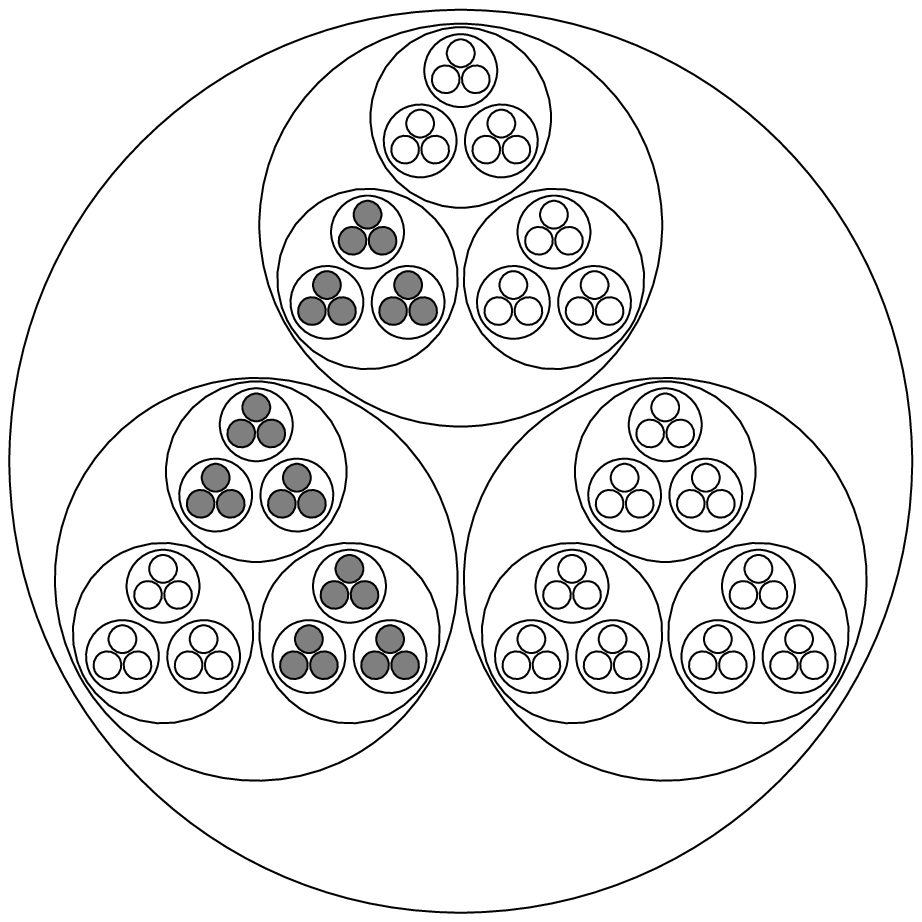}}
\caption{$\Omega_{1,0}$ and $\Omega_{1,1}$ of Example~\ref{ssec:Q3wave}.}
\label{fig:Omega5.3b}
\end{figure}

\begin{figure}
\scalebox{.8}{\includegraphics{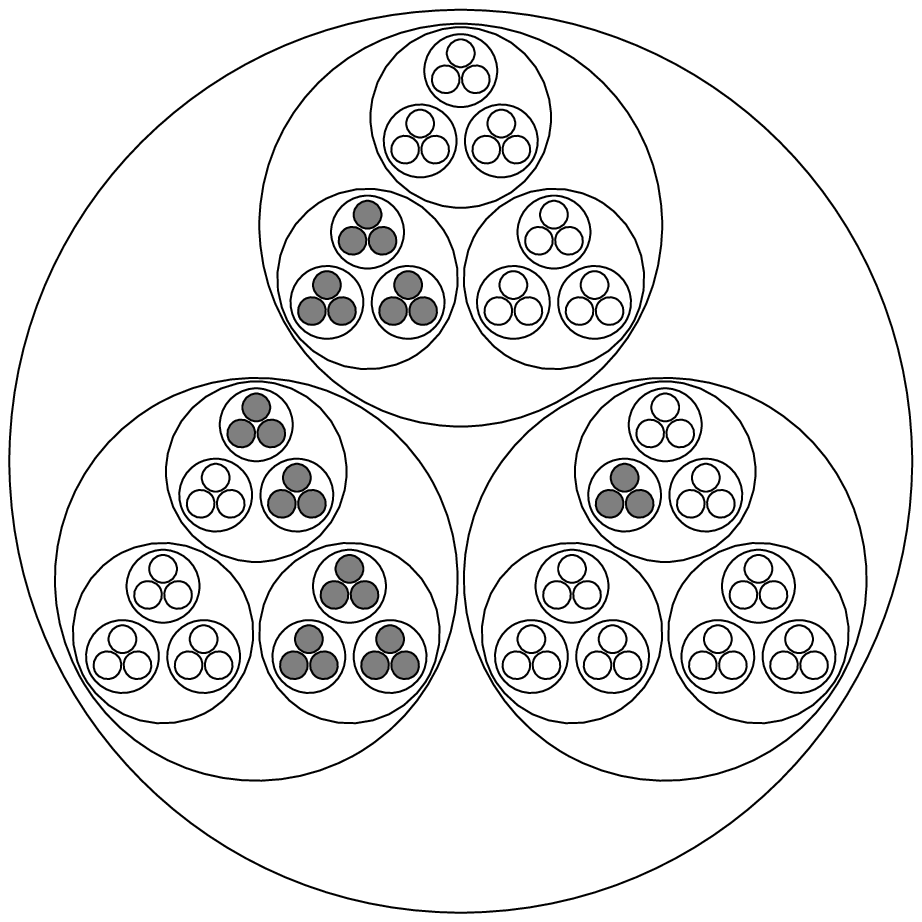}
\includegraphics{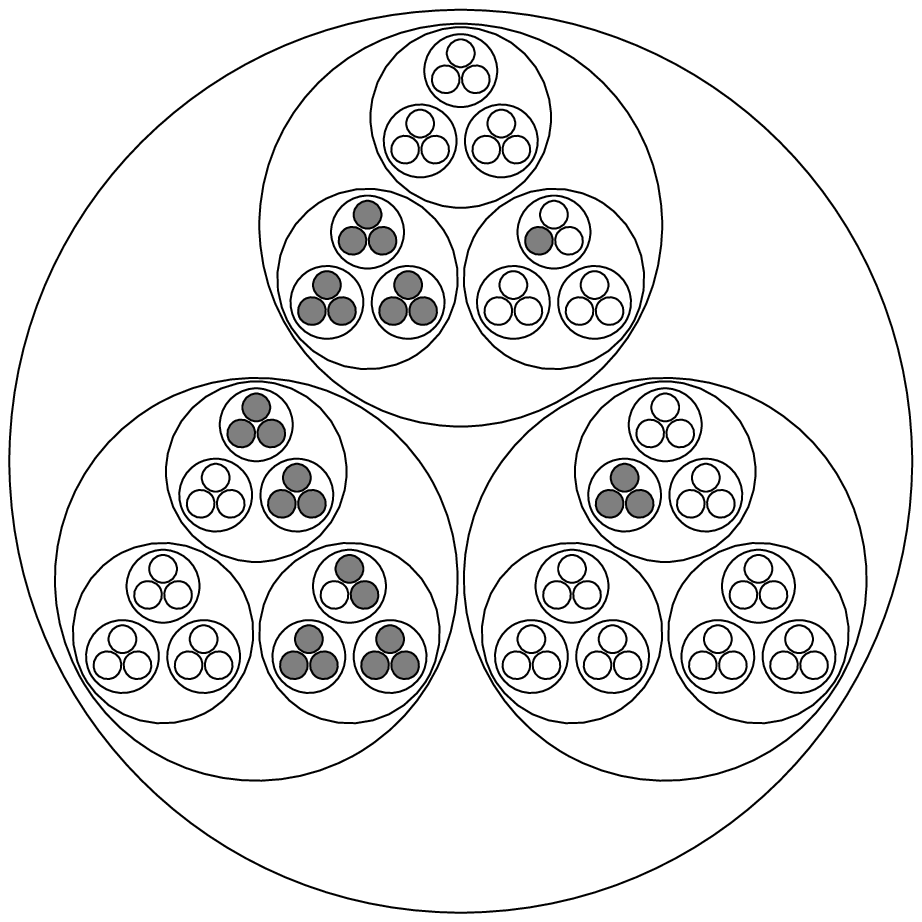}}
\caption{$\Omega_{1,2}$ and $\Omega_{1,3}$ of Example~\ref{ssec:Q3wave}.}
\label{fig:Omega5.3e}
\end{figure}

Take $M=0$, so that $W=H^{\perp}$, set $N=1$,
and let $\Omega_{1,0}=H^{\perp}$.
For $\gamma\in H^{\perp}$, define
$$
T_1(\gamma) =
\begin{cases}
\gamma + 2/3 & \text{if } \gamma\in 1+3\ZZ_3, \\
\gamma + 1/3 & \text{if } \gamma\in (3\ZZ_3) \cup (2 + 3\ZZ_3),
\end{cases}
$$
as in Figure~\ref{fig:Tmap5.3}.
Clearly $\Lambda_{1,1}=3\ZZ_3$, which is a disk of measure $1/3$.
By induction, it is easy to show that for $m\geq 1$,
$$
\Lambda_{1,m} =
\begin{cases}
  2 +1\cdot 3 + 2\cdot 3^2 + 1\cdot 3^3 + \cdots
  + 1\cdot 3^{m-2} + 3^m\ZZ_3
  &\text{if $m$ is odd}; \\
  1 +2\cdot 3 + 1\cdot 3^2 + 2\cdot 3^3 + \cdots
  + 1\cdot 3^{m-2} + 3^m\ZZ_3
  & \text{if $m$ is even}.
\end{cases}
$$
Note that the sums stop at $3^{m-2}$, so that
$\Lambda_{1,1}=3\ZZ_3$ and $\Lambda_{1,2}=1 + 3^2 \ZZ_3$.
Equivalently,
$$
\Lambda_{1,m} =
\begin{cases}
  -5/8 + 3^{m-1} + 3^m\ZZ_3
  &\text{if $m$ is odd}; \\
  -7/8 + 3^{m-1} + 3^m \ZZ_3,
  & \text{if $m$ is even}.
\end{cases}
$$
As in the previous example,
$\Lambda_{1,m}$ is a disk of measure $1/3^{m}$,
but it lies in $3\ZZ_3$ if $m=1$,
in $1+3\ZZ_3$ if $m$ is even,
and in $2+3\ZZ_3$ if $m\geq 3$ is odd.
See Figures~\ref{fig:Omega5.3b}--\ref{fig:Omega5.3e}
for some of the resulting sets $\Omega_{1,m}$.
Thus, 
$\Lambda_1=\bigcup_{m\geq 1} \Lambda_{1,m} \subseteq\ZZ_3$
is a disjoint union of countably many
disks, with $\nu(\Lambda_1) = 1/(3-1)=1/2$, and
the wavelet set is
\begin{multline*}
\Omega_1 =
\left[ \ZZ_3 \setminus \bigcup_{n=1}^\infty 
\left( (-5/8 + 3^{2n-2} + 3^{2n-1} \ZZ_3) \cup
(-7/8 + 3^{2n-1} + 3^{2n} \ZZ_3) \right) \right] \\
{} \cup
\bigcup_{m=1}^\infty 
\left( (-7/24 + 3^{2n-2} + 3^{2n-1} \ZZ_3) \cup
(-5/24 + 3^{2n-1} + 3^{2n} \ZZ_3) \right).
\end{multline*}

ACKNOWLEDGEMENTS:
Many of the references reflect
fruitful conversations about our program with the authors, both analysts and
number theorists. In particular, and not necessarily
referenced, we want to thank Professors Lawrence Baggett,
David Cox, Hans Feichtinger, Karlheinz Gr{\"o}chenig, and
Niranjan Ramachandran.

\nocite{*}

\bibliographystyle{amsplain}
\bibliography{BandB}


\end{document}